\input amstex.tex
\documentstyle{amsppt}  
\input graphicx.tex 

\magnification=1200
\hsize=150truemm
\vsize=224.4truemm
\hoffset=4.8truemm
\voffset=12truemm

\TagsOnRight
\NoBlackBoxes
\NoRunningHeads

\def\Square{\rlap{$\sqcup$}$\sqcap$}
\def\cqfd {\quad \hglue 7pt\par\vskip-\baselineskip\vskip-\parskip
{\rightline{\Square}}}

\define\E{{\Cal E}}
\redefine\P{{\Cal P}}
\define\T{{\Cal T}}
\redefine\D{{\Cal D}}
\define\G{{\Cal G}}
\define\Q{{\bold Q}}
\define\R{{\bold R}}
\define\Z{{\bold Z}}
\define\mi{^{-1}}
\let\thm\proclaim
\let\fthm\endproclaim
\let\inc\subset 
\let\ds\displaystyle
\let\ov\overline
\define\Aut{\text{\rm Aut}}
\define\Out{\text{\rm Out}}

\newcount\tagno
\newcount\secno
\newcount\subsecno
\newcount\stno
\global\subsecno=1
\global\tagno=0
\define\ntag{\global\advance\tagno by 1\tag{\the\tagno}}

\newcount\figno
\newcount\fihno
\global\figno=0
\global\fihno=1
\define\fig{\global\advance\figno by 1 \global\advance\fihno by 1
{\the\figno}}
\define\sta{\the\secno.\the\stno
\global\advance\stno by 1}

\define\sect{\global\advance\secno by
1\global\subsecno=1\global\stno=1\
\the\secno. }

\def\nom#1{\edef#1{\the\secno.\the\stno}}
\def\eqnom#1{\edef#1{(\the\tagno)}}

\newcount\refno
\global\refno=0
\def\nextref#1{\global\advance\refno by 1\xdef#1{\the\refno}}
\def\bref {\ref\global\advance\refno by 1\key{\the\refno}}

\nextref\BK
\nextref\BS
\nextref\Cld
\nextref\Clc
\nextref\Coll
\nextref\Collev
\nextref\CV
\nextref\Fe
\nextref\Fed
\nextref\Forr
\nextref\Foc
\nextref\For
\nextref\GHMR
\nextref\Gu
\nextref\GL
\nextref\GLL
\nextref\Krn
\nextref\Kr
\nextref\KV
\nextref\LGD
\nextref\Ler
\nextref\Lrang
\nextref\LL
\nextref\McM
\nextref\Mold
\nextref\MSW
\nextref\Pett
\nextref\Pet 
\nextref\RV
\nextref\RR
\nextref\Wh

\topmatter

\abstract   A generalized Baumslag-Solitar group (GBS group) is a finitely
generated group $G$ which acts on a tree with all edge and vertex stabilizers
infinite cyclic. We show that Out$(G)$ either contains non-abelian
free groups or is virtually nilpotent of class $\le2$. It has torsion only at
finitely many primes.

One may decide algorithmically whether  Out$(G)$ is virtually nilpotent or
not. If it is, one may decide whether  it is virtually abelian, or
finitely generated. The isomorphism problem is solvable among GBS groups
with Out$(G)$ virtually nilpotent.

If $G$ is unimodular (virtually $F_n\times{\bold Z}$), then Out$(G)$ is
commensurable with a semi-direct product ${\bold Z}^k\rtimes\text{\rm
Out}(H)$ with
$H$ virtually free.   

\endabstract

\title On the automorphism group of generalized Baumslag-Solitar groups
 \endtitle

\author  Gilbert Levitt
\endauthor

\toc
\widestnumber\head{0}
\head 1. Introduction and statement of results\endhead
\head 2. Basic facts about GBS groups\endhead
\head 3. The automorphism group of a GBS tree\endhead
\head 4. Unimodular groups\endhead
\head 5. The deformation space\endhead
\head 6. Free subgroups in $\Out(G)$\endhead
\head 7. Groups with $\Out(G)\not\supset F_2$\endhead
\head 8. Further results\endhead
\head  { } References\endhead
\endtoc

\endtopmatter 

\document

\head \sect Introduction and statement of results  \endhead

The groups $BS(m,n)=\langle a,t\mid ta^mt\mi=a^n\rangle$ were introduced
by Baumslag-Solitar   [\BS] as very simple examples of
non-Hopfian groups (a group $G$ is non-Hopfian if there exists a non-injective
epimorphism from $G$ to itself). It is now known that 
$BS(m,n)$ is Hopfian if and only if $m=\pm1$, or $n=\pm1$, or $m,n$ have the
same set of prime divisors [\Collev]. In particular,
$BS(2,4)$ is Hopfian while $BS(2,3)$ is not.

Though it has exotic epimorphisms, $BS(2,3)$ has very few automorphisms:
its automorphism group is generated by inner automorphisms and
the obvious involution sending $a$ to $a\mi$ [\Coll, \GHMR]. On the other hand,
$BS(2,4)$ has an incredible number of automorphisms, as {\it its
automorphism group is not finitely generated\/} [\Collev].

The reason behind this drastic difference is that, because 2 divides 4 but not
3, the presentation of
$BS(2,4)$ is much more flexible than that of $BS(2,3)$. By this we mean, in
particular, that $BS(2,4)$ admits the   infinite sequence of
presentations
$$BS(2,4)=\langle   a, b,t\mid tb t\mi=b{}^2,
b^{2^p}=a^2\rangle\tag{$1_p$}$$  
obtained from the standard one by introducing a new generator
$b=t^{-p}a^2t^p$. It is  clear already from $(1_p)$ that
$G=BS(2,4)$ has many automorphisms, as fixing $b ,t$ and
conjugating $a$ by $b$  defines an element of order $2^p$ in
$\Out(G)$.

The presentations $(1_p)$ express $BS(2,4)$ as a generalized
Baumslag-Solitar group, or {\it GBS group\/}, or graph of $\Z$'s, namely as
the  fundamental group of a finite graph of groups $\Gamma $ with all edge and
vertex groups infinite cyclic. This is visualized as a {\it labelled graph\/}, with
the absolute value of the labels indicating the index of edge groups in vertex
groups (see Figure \fig).

\midinsert
\centerline 
{\includegraphics[scale=.5]
{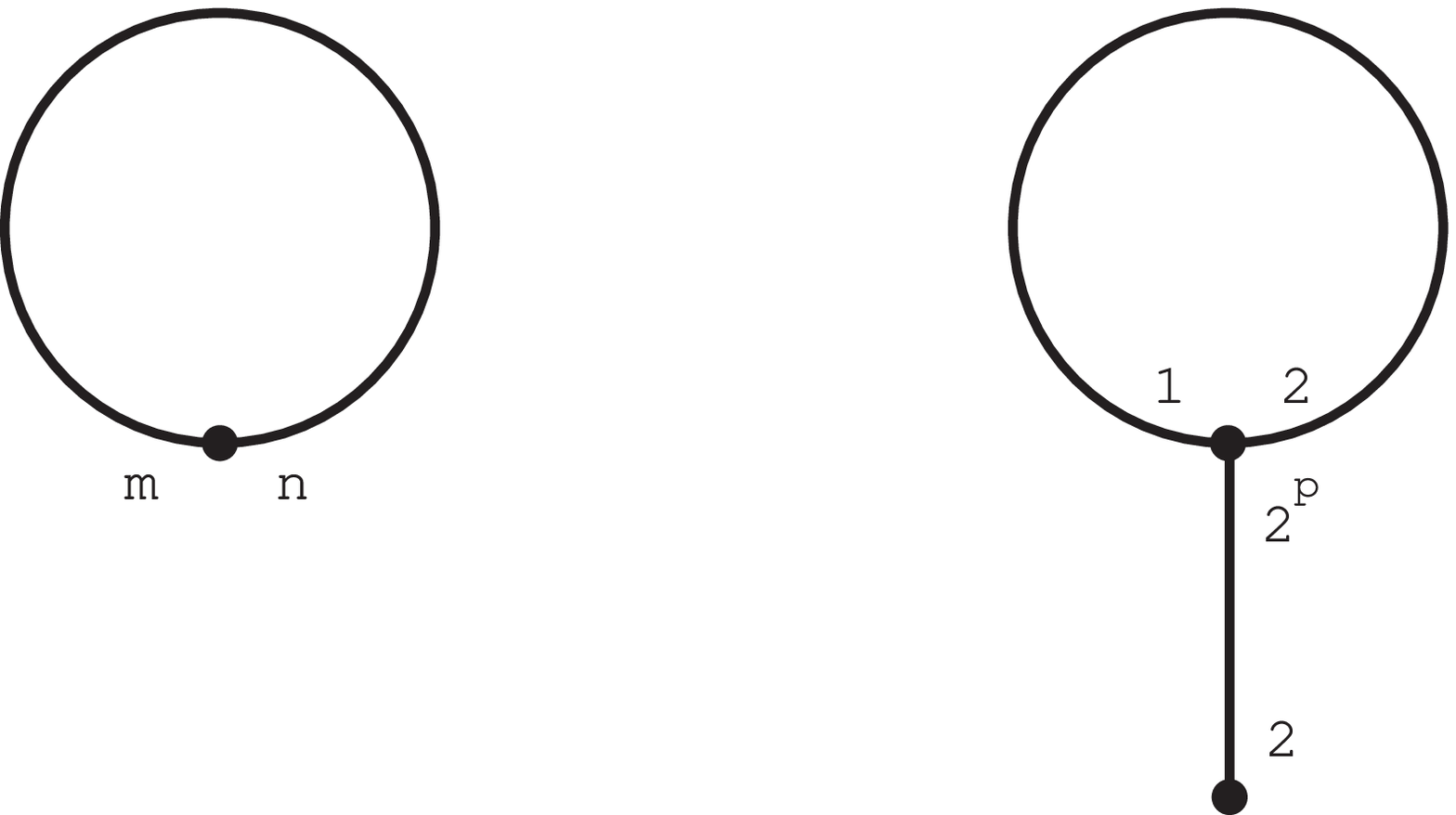}}
\captionwidth{220pt}
\botcaption 
 {Figure \the\figno}
{The labelled graphs associated to the standard presentation of
$BS(m,n)$, and to $(1_p)$. }
\endcaption
\endinsert

In this paper, we study automorphisms of   GBS  groups. See [\Foc, \For, \Kr,
\Lrang,
\Wh]  for various algebraic and geometric properties of these groups. As
pointed out in [\Foc], they are especially interesting in connection with  JSJ
theory. 

Before giving general results,
let us review certain classes of GBS groups for which more specific
statements may be obtained. They are defined either by ``local'' conditions on
the labelled graph, or by ``global'' algebraic conditions on the group. In the
rest of this introduction, we always assume that $G$ is not one of the
elementary GBS groups: $\Z$,
$\Z^2$, and the Klein bottle group. 

\subhead{Algebraically rigid groups}\endsubhead

As evidenced by the example of $BS(2,4)$, the main difficulty with GBS groups
is that they may be represented by  many different labelled
graphs $\Gamma $.  Sometimes, though, $\Gamma $ is essentially unique. By
[\GHMR, \Pet ], this {\it algebraic rigidity\/}  holds  in particular  when there is
no {\it divisibility relation\/} in
$\Gamma
$: if $p,q$ are labels near the same vertex, then $p$ does not divide $q$
(see Section 2 for a precise definition and a characterization of algebraic
rigidity). 

Given $\Gamma $, let $T$ be the associated  Bass-Serre tree, which we call a
{\it GBS tree\/}. Let $\Out^T(G)\inc \Out(G)$ be
the  subgroup leaving $T$ invariant. Most elements of $\Out^T(G)$
may be viewed as ``twists'' (see Section 3). Algebraic rigidity implies
$\Out^T(G)=\Out(G)$, but in general $\Out^T(G)$ is smaller. 

\nom\un
\thm{Theorem \sta} Let $G$ be a GBS group, represented by a labelled
graph $\Gamma $, and let
$T$ be the Bass-Serre tree. Define $k$ as the first Betti number $b$ of
$\Gamma
$ if $G$ has a nontrivial center, as $b-1$ if the center is trivial.
\roster
\item  The
torsion-free rank of the abelianization of $G$ is $k+1$. 
\item
The group $\Out^T(G)$ is virtually $\Z^k$.
\item Up to
commensurability within $\Out(G)$, the subgroup  $\Out^T(G)$  
does not depend on
$\Gamma $.
\endroster
\fthm

Conversely, any subgroup of $\Out(G)$
commensurable with  a subgroup of
$\Out^T(G)$ is contained in $\Out^{T'}(G)$ for some GBS tree $T'$ [\Clc].

For $G=BS(m,n)$, one has $k=0$ if $m\neq n$,
and
$k=1$ if $m=n$. For $G=BS(2,4)$ with the presentation $(1_p)$, the group
$\Out^T(G)$ has order $2^{p+1}$. 

\thm{Corollary \sta}  If $G$ is algebraically rigid, then $\Out(G)$ is virtually
$\Z^k$. 
\fthm

The converse is also true (see Theorem 8.5).

\subhead Unimodular groups\endsubhead

A GBS group $G$ is {\it unimodular\/} if $xy^px\mi=y^q$ with $y\neq1$ implies
$|p|=|q|$, or equivalently if $G$ is virtually $F_n\times\Z$ (with $F_n$ a free
group of rank $n$). The group $G$ then has a normal infinite cyclic subgroup
with virtually free quotient, and we show: 

\thm{Theorem \sta} If $G$ is  unimodular, there  is a split
exact sequence
 $$\{1\}\to
\Z^{k}\to\Out_0(G) \to \Out_0(H)
\to\{1\},$$ where $k$ is as above, $H$ is virtually free, and 
 $\Out_0$ has finite index in $\Out$.
\fthm

Since $\Out(H)$ is VF [\KV], we get:

\thm{Corollary \sta} $\Out(G)$ is virtually torsion-free and VF (it has a finite
index subgroup admitting a finite classifying space).
\fthm
 
\subhead Groups  with no nontrivial integral modulus\endsubhead

Now consider groups $G$ which do not contain a solvable Baumslag-Solitar
group
$BS(1,n)$ with $n\ge2$ (there is an equivalent characterization   in terms of
the modular homomorphism  $\Delta :G\to\Q^*$, see Section 2). 

Given any GBS group $G$, the group $\Out(G)$ acts on the space $P\D$ of all
GBS trees (see Section 5), with stabilizers virtually $\Z^k$ by Theorem \un.
Clay [\Cld] proved that the   space $P\D$ is contractible (see also [\GLL]), and
Forester [\For] proved that the quotient is a finite complex if $G$ does not
contain 
$BS(1,n)$ with $n\ge2$.  This gives:

\thm{Theorem \sta} If $G$ does not contain 
$BS(1,n)$ for $n\ge2$, then $\Out(G)$
is $F_\infty$ (in particular, it is finitely presented). If furthermore $\Out(G)$
is virtually torsion-free, then it  is VF.
\fthm

\subhead{Arbitrary groups}
\endsubhead

Now let $G$ be any GBS group. 

\nom\six
\thm{Theorem \sta} Either $\Out(G)$ contains a nonabelian free group, or it is
virtually nilpotent of class $\le2$. 
\fthm

A group is nilpotent of class $\le2$ if and only if every commutator is central. 
As an example, let $G=\langle a,s,t\mid sa=as, ta^2=a^2t\rangle$ (see Figure
\fig). Then
$\Out(G)$ is virtually the integral Heisenberg group $H_3$, with 
$\pmatrix1&i&j\\0&1&k\\0&0&1\endpmatrix$ mapping $(a,s,t)$ to
$(a,sa^k,ts^ia^j)$. 

\midinsert
\centerline 
{\includegraphics[scale=.5]
{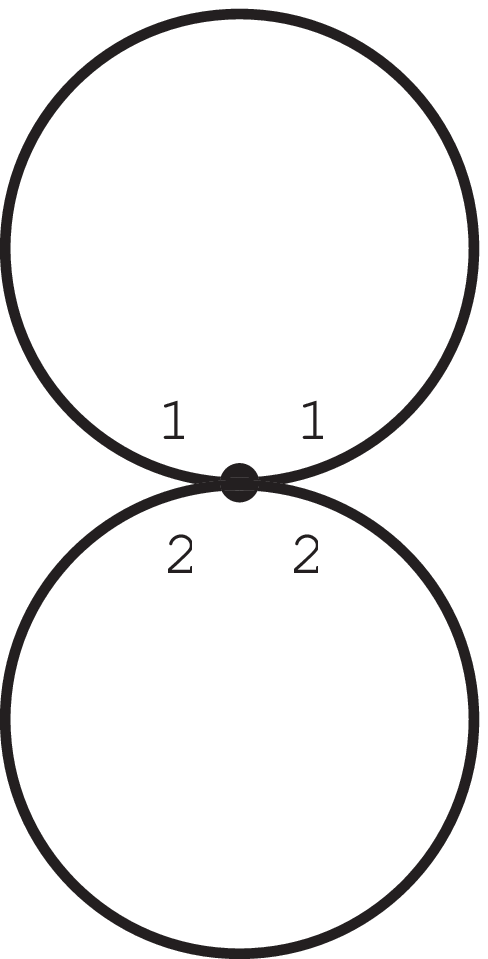}}
\botcaption
 {Figure \the\figno}
{$\Out(G)$ is virtually the integral Heisenberg group $H_3$. }
\endcaption
\endinsert

Which possibility of Theorem \six{} occurs may be explicitly decided from 
the divisibility relations in
  any labelled graph $\Gamma $ representing $G$. We have seen that $\Out(G)$
is virtually abelian if there is none. A key observation is that certain
divisibility relations force the existence of $F_2$ inside $\Out(G)$. 

\midinsert
\centerline 
{\includegraphics[scale=.5]
{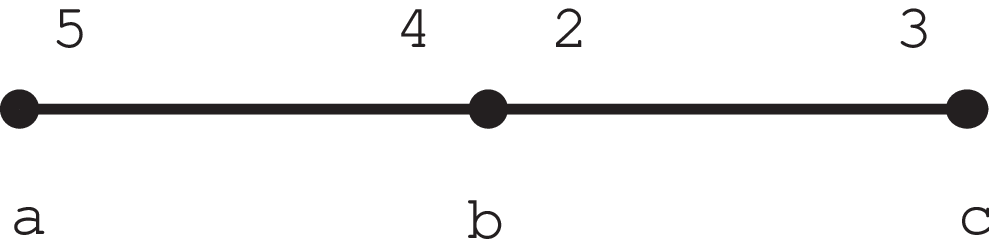}}
\captionwidth{220pt}
\botcaption
 {Figure \the\fihno}
{$\Out(G)$ contains $F_2$. }
\endcaption
\endinsert

As a basic example, consider $G=\langle
a,b,c\mid a^5=b^4, b^2=c^3\rangle$ (see Figure
\fig). It is the amalgam of $G_1=\langle
a\rangle$ with $G_2=\langle b,c\rangle$ over $C=\langle b^4\rangle$. The
divisibility relation $2|4$ at the middle vertex implies that $C$ is central in
$G_2$. For any
$g\in G_2$, we may therefore define an automorphism $\varphi _g$ of $G$ as
being the identity on
$G_1$ and conjugation by $g$ on $G_2$. It is easy to show that the subgroup of
$\Out(G)$ generated by the $\varphi _g$'s is isomorphic to $\langle
b,c\mid   b^2=c^3=1\rangle$, hence contains $F_2$.

To prove Theorem \six{}, we assume that $\Out(G)$ does not contain $F_2$ and
 we   describe   which divisibility relations may occur (Section 6).
 In Section 7, we show that, though the
GBS tree
$T$ may not be
$\Out(G)$-invariant, some (non GBS) tree $S$ obtained from $T$ by collapsing
certain edges is. We then prove  that $\Out^S (G) $
is virtually nilpotent. 

\midinsert
\centerline 
{\includegraphics[scale=.5]
{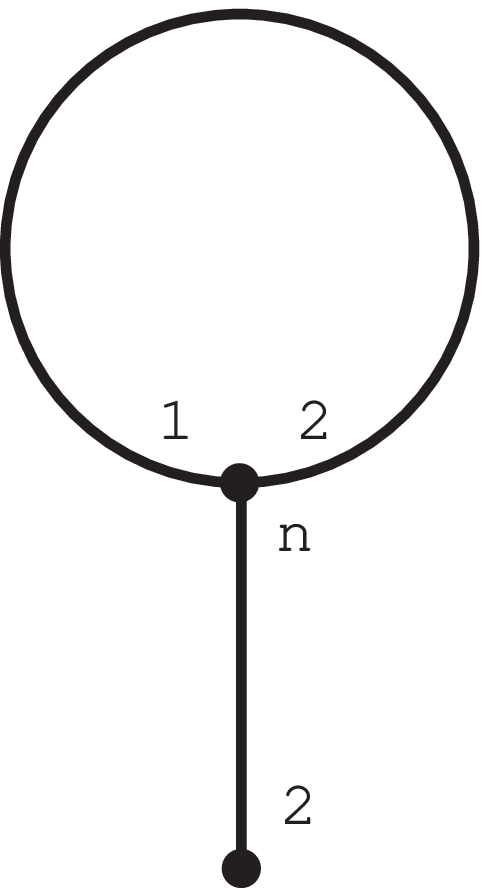}} 
\captionwidth{210pt}
\botcaption
 {Figure \the\fihno}
{$\Out(G)$ is virtually an infinitely generated abelian group if $n$ is not a
power of
$2$. }
\endcaption
\endinsert

For instance, let $G=\langle   a,b,t\mid tbt\mi=b^2,b^n=a^2\rangle$ with $n$
not a power of 2 (see Figure
\fig). In this case,
$S$ is obtained from $T$ by collapsing edges projecting onto the loop of
$\Gamma
$. The group
$\Out(G)$ is virtually abelian (but not finitely generated).

A special case of Theorem \six{} is:

\thm{Theorem \sta} If no label of $\Gamma $ equals 1, then $\Out(G)$ contains
$F_2$ or is a finitely generated   virtually abelian group.  
\fthm

As a corollary of our analysis, we show:

\thm{Theorem \sta}  The isomorphism problem is solvable for GBS groups $G$
such that $\Out(G)$ does not contain $F_2$.
\fthm

The isomorphism problem for GBS groups is to decide whether two
labelled graphs define isomorphic groups. It is solvable for groups with no
nontrivial integral modulus [\For] and $2$-generated groups [\Lrang], but
open in general.

We also show:

\nom\huit
\thm{Theorem \sta} The set of prime numbers $p$ such that $\Out(G)$
contains non-trivial $p$-torsion is finite. 
\fthm

The paper is organized as follows. In Section 2, we review basic
properties of GBS groups, such as algebraic rigidity and the modular
homomorphism $\Delta  $. We extend to GBS groups a 
result of [\Fed] about twisted conjugacy classes.
In Section 3, we study $\Out^T(G)$, proving Theorems \un{}  and \huit. Section
4 is devoted to unimodular groups, Section 5 to the action of $\Out(G)$ on
$P\D$. Theorem \six{} is proved in Sections  6 and 7.    In Section 8 we discuss
several special cases and the isomorphism problem.  

\smallskip

\vskip8pt
{\baselineskip=9pt
{\eightpoint Acknowledgments. {\it I am grateful to M\. Forester and M\. Clay
for helpful suggestions, and to P\. Papasoglu for pointing out that
$\Out(BS(2,4))$ is not finitely generated.}}\par} 

\head \sect Basic facts about GBS groups  \endhead

\subhead   Labelled graphs  \endsubhead

A GBS group $G$ is the fundamental group of a finite graph of groups
$\Gamma $ whose vertex and edge groups are all infinite cyclic. It is
torsion-free. 

We denote by $b$ the first Betti number of the graph $\Gamma $. Note
the distinction between $G= \pi _1(\Gamma )$ and the topological
fundamental group $\pi _1^{top}(\Gamma )\simeq F_b$. 

If we choose generators for edge and vertex groups, the inclusion
maps are multiplications by non-zero integers. An oriented edge $e$
thus has a label $\lambda _e\in\Z\setminus\{0\}$, describing the inclusion of
the edge group $G_e$ into the vertex group $G_{o(e)}$  at the origin of $e$.
As in [\For], we visualize the graph of groups as a {\it labelled graph\/}
$\Gamma $, with the label $\lambda _e $  pictured near the origin
$o(e)$. 

A pair
$\varepsilon =(e,\ov e)$ of opposite edges is a non-oriented edge. It  
carries two labels, one near either endpoint  $o(e),o(\ov e)$, and we
say that $\varepsilon $ (or $e$) is a {\it $(\lambda _e,\lambda _{\ov
e})$-edge\/}.  An edge is a   {\it loop\/} if its endpoints are
equal, a {\it segment\/} if they are distinct. 

The group $G$ associated to a labelled graph $\Gamma $ may be
presented as follows.   Choose a maximal subtree $\Gamma _0\inc
\Gamma $.  There is one generator $x_v$  for each vertex $v$, and one
generator $t_\varepsilon $ for each non-oriented edge
$\varepsilon$ not in $\Gamma _0$. Each   non-oriented edge
$\varepsilon =(e,\ov e)$ of $\Gamma $
  contributes one relation. If $ \varepsilon $ is contained in
$\Gamma _0$, the relation is $(x_{o(e)})^{{\lambda (e)}}=  (x_{o(\ov
e)})^{{\lambda (\ov e)}}$. If $\varepsilon  $ is not in $\Gamma _0$, the
relation is
$t_e(x_{o(e)})^{{\lambda (e)}}t_e\mi=   (x_{o(\ov e)})^{{\lambda (\ov
e)}}$.

Replacing     the    chosen  generator of a vertex group $G_v$ by its
inverse   changes the sign of   all labels near  $v$.  Replacing an edge
group generator   changes the sign of   both  labels   carried by the
edge. These changes are {\it admissible sign changes\/}.
Labelled graphs will  always be considered up to admissible sign
changes. 

When we  focus on a particular edge, we always use admissible sign changes
to make it a $(p,q)$-segment with $p,q>0$, or a $(p,q)$-loop
with $1\le p\le|q|$.  

A
$(1,q)$-loop is an {\it ascending loop\/}. It is  a {\it strict ascending
loop\/} if $|q|>1$; note that $G$ then contains a solvable
Baumslag-Solitar group
$BS( 1,  q)$. A
$(p,q)$-loop with $p|q$ is a {\it pseudo-ascending loop\/}.

\subhead   GBS trees \endsubhead

 Let $G$ be the fundamental group of a labelled graph $\Gamma $. The
associated Bass-Serre tree is a locally finite $G$-tree $T$ with all
edge and vertex stabilizers infinite cyclic. Such $G$-trees will be
called {\it GBS trees\/}.  Two trees are considered to be the same if
there is a $G$-equivariant isomorphism between them.

 We always assume that the action is {\it minimal\/}: there is no proper
invariant subtree. In terms of $\Gamma $, this is equivalent to saying
that the label  near every terminal vertex is bigger than $1$. We also
assume that actions are without inversions.  

Given a   GBS tree $T$, one obtains a labelled graph $\Gamma
=T/G$, with the labelling   well-defined up to admissible sign
changes (see Remark 2.3 in [\For]). This graph of groups is {\it marked\/}:
there is an isomorphism from its fundamental group to $G$, well-defined up to
composition with an inner automorphism.  The valence of a vertex   $v\in T$  
is the sum of the absolute values of the labels near its image in $\Gamma $.

GBS trees $T$ and marked labelled graphs $\Gamma $ are thus
equivalent concepts. We will work with both. We usually use the same letter
$v$ (resp\.
$e$) for a vertex (resp\.
edge) of $T$ and its image in $\Gamma $. When we need to distinguish, we
write $\ov v$ for the image of $v$ in $\Gamma $. 
We denote vertex stabilizers
(vertex groups) by $G_v$, edge stabilizers (edge groups) by $G_e$.

\subhead   Collapses and algebraic
rigidity\endsubhead

Collapsing an edge $e$ of $\Gamma $ (or equivalently a $G$-orbit of edges of
$T$) yields a new tree $S$, which usually is not a GBS tree. It is a GBS tree if
and only if $e$ is a segment and at least one of the labels $\lambda _e$,
$\lambda _{\ov e}$ equals $1$. Such an edge will
be called a {\it collapsible\/} edge. 

In the proof of Theorem \six, we will  
collapse
$(2,2)$-edges and $(1,q)$-loops; these are not collapsible edges. We usually
denote by $\Theta $ the collapsed graph of groups, by $\pi :T\to S$ the
collapse map.  The image of a vertex $v\in T$ is denoted by $\pi
(v)$, or sometimes just
$v$. The  stabilizer of $\pi (v)$ in $S$ contains the cyclic group $G_v$, we call
it 
$H_v$ (it will often be a solvable Baumslag-Solitar group). We use the same
letter for a non-collapsed edge   of
$T$ and its image in $S$. It has the same stabilizer in both trees.

Collapsing a collapsible edge  is called an {\it elementary   collapse\/}. The
reverse move is an {\it elementary  expansion\/}.
Labels near
$o(  e)$ get multiplied by
$\lambda _{\ov e}  $ when we collapse an edge $e$ with $\lambda _e=1$ (see
Figure
\fig). 

\midinsert
\centerline 
{\includegraphics[scale=.5]
{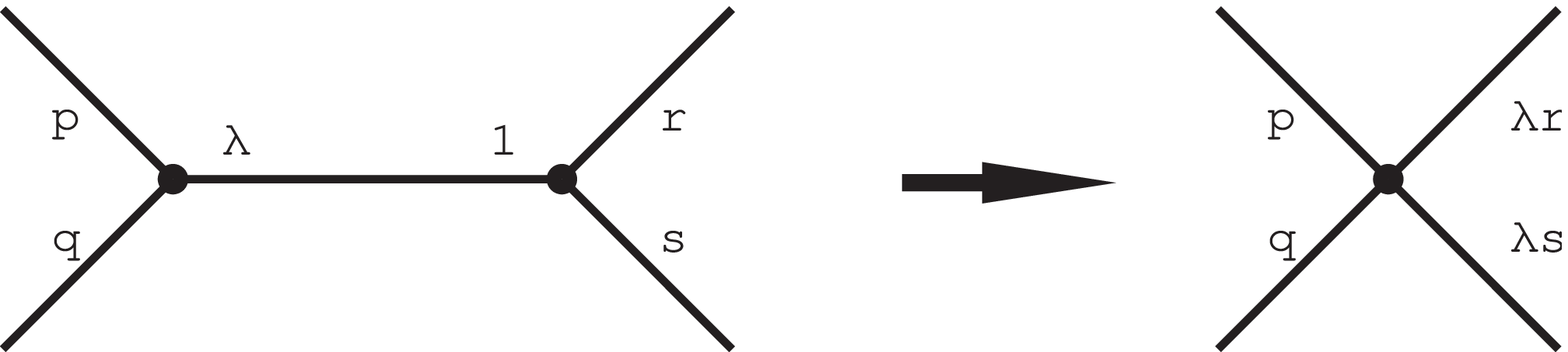}}
\captionwidth{220pt}
\botcaption
 {Figure \the\figno}
{Elementary collapse. }
\endcaption
\endinsert

The graph $\Gamma  $, or the tree $T$,  is {\it reduced\/} (in the sense
of [\Forr]) if there is no collapsible edge. 
In terms of trees, $T$ is reduced if and only if
any edge $e=vw$ satisfying $G_e=G_v$  has its endpoints in the same
$G$-orbit.  Any tree may be reduced by applying a finite
sequence of elementary collapses (the reduction is not always unique).

\midinsert
\centerline 
{\includegraphics[scale=.5] 
{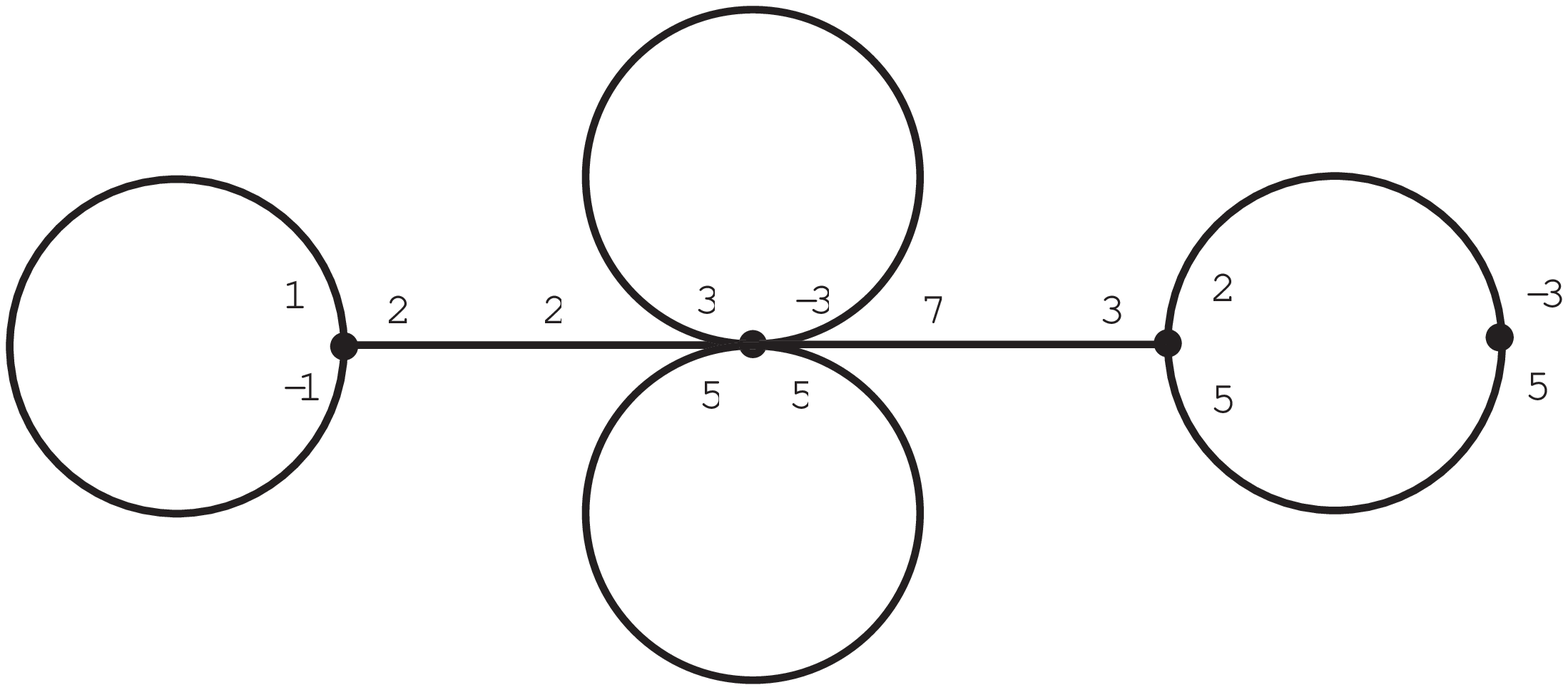}}
\botcaption
 {Figure \the\fihno}
{A labelled graph representing an algebraically rigid group.}
\endcaption
\endinsert

A reduced GBS tree $T$ is {\it rigid} if it is the only reduced GBS
 tree (up to equivariant isomorphism). Building on  work from [\Forr,
\GHMR, \Gu,
\Pet], it is shown in [\Ler] that, if $G$ is not solvable,  $T$ is rigid if and only if
$\Gamma
$ satisfies the following condition (see Figure
\fig): if $e,f$
are distinct oriented edges of $\Gamma $ with the same origin
$v$, and the label of $f$ near $v$ divides that of $e$, then either $e=\ov f$ is a
$(p,\pm p)$-loop with $p\ge2$, or $v$ has valence 3 and bounds a
$(1,\pm1)$-loop.

  In particular,
$T$ is rigid  whenever  there is no divisibility relation in $\Gamma $ (recall
that a  {\it divisibility relation\/}  is a relation $p|q$ between two labels at the
same vertex). 

When there is a rigid GBS tree, we say that $G$ is {\it algebraically rigid\/}. In
this case, there is only one reduced marked labelled graph  
representing
$G$. See [\MSW, \Wh] for quasi-isometric rigidity of GBS groups.

\subhead   Non-elementary groups \endsubhead

We say that $G$ is {\it elementary\/} if
$T$ may be chosen to be a point or a line. As vertices of $T$ then have
valence at most 2, there are only four possibilities for $\Gamma $: a
point, a
  $(1,1)$-loop, a
  $(1,-1)$-loop, a $(2,2)$-segment. The corresponding groups are
 $\Z$,   $\Z^2$, and   the Klein bottle group $\langle x,t\mid
txt\mi=x\mi\rangle=\langle a,b\mid a^2=b^2\rangle$, with $\Out(G)$ equal to
$\Z/2\Z$, $GL(2,\Z)$, and $\Z/2\Z\times\Z/2\Z$ respectively.

Though non-elementary, the solvable groups $BS(1,n)$ are special. For $|n|>1$,
the group $\Out(BS(1,n))$ is virtually $\Z^{r-1}$, where $r$ is the number of
prime divisors of $n$ [\Coll]. More generally, see [\Coll, \Collev, \GHMR] for a
presentation of $\Out(BS(m,n))$.

From now on, we consider only non-elementary groups. Here are a few
simple properties (compare [\Foc]).

From the action of $G$ on $T$, it is easy to see that a non-elementary
GBS group either is a solvable Baumslag-Solitar group $BS(1,n)$ (if it
fixes an end of
$T$), or contains non-abelian free groups (if the action on $T$ is
irreducible). In particular, $G$ always has exponential growth. A finitely
generated subgroup of
$G$ is free (if it acts freely on $T$) or is a GBS group. A non-elementary
GBS group is  one-ended, coherent, and has cohomological dimension 2 [\Foc,
\Kr].

Any GBS group maps onto $\Z$ 
(the   presentation of $G$ given earlier has more generators than
relators).   GBS groups are therefore locally indicable, hence  orderable (see
[\RR]).  Using the fact that $\langle a,b\mid a^p=b^q
\rangle$ is bi-orderable only if $|p|$ or $|q|$ equals $1$, one shows that the
only   bi-orderable GBS groups are  
$F_n\times\Z$ and  $BS(1,n)$ for $n\ge1$.

\subhead   Elliptic elements \endsubhead

Two subgroups $H,K$ of $G$ are {\it commensurable\/} if $H\cap K$ has
finite index in both $H$ and $K$. Two elements $g,h$ are commensurable
if $\langle g
\rangle$ and $\langle h
\rangle$ are commensurable, equivalently if there is a relation
$g^p=h^q$ with
$p,q$   non-zero integers. The {\it commensurator\/} of $g$ is the
subgroup
$Comm(g)$ consisting of all $x\in G$ such that $xgx\mi$ is
commensurable to $g$. 

Given any $G$-tree $T$, an element $g\in G$, or a subgroup $H$, is {\it
elliptic\/}   if it fixes a point.  If
$g$ is not elliptic, it is {\it hyperbolic\/}: there is an invariant axis, on
which $g$ acts as a translation by some positive integer $\ell(g) $.
Conjugate or commensurable elements have the same type (elliptic or
hyperbolic). A relation
$ga^pg\mi=a^q$ with
$|p|\ne|q|$ implies that $a$ is elliptic, because its translation length
satisfies 
$|p|\ell(a)=|q|\ell(a)$.

\thm{Lemma \sta{} [\Forr]} Let $T$ be a   GBS  tree, with $G$
non-elementary.  Any two  non-trivial elliptic elements
$g,g'$ are commensurable. An element $g\in G$ is elliptic if and only if
its commensurator equals $G$.
\fthm 

\demo{Proof} If $g,g'$ fix vertices $v,v'$, one shows that they are
commensurable by induction on the distance between $v$ and $v'$. If
$g$ is hyperbolic, its axis is $Comm(g)$-invariant, so $Comm(g)\neq G$
because $G$ is not elementary.
\cqfd\enddemo

\nom\ellind
\thm{Corollary \sta}  The set of elliptic elements depends only on $G$,
not on the GBS tree $T$. It is invariant under automorphisms of $G$.
\cqfd
\fthm

As any two GBS trees have the same elliptic subgroups, Forester's
deformation theorem  [\Forr] yields: 

\nom\defor
\thm{ Corollary \sta {}} Let $G$ be a non-elementary GBS group.
Any two GBS
 trees are related (among GBS trees) by a finite sequence of elementary
expansions and collapses. \cqfd
\fthm

The quotient $G/G_{ell}$ of $G$ by the subgroup generated by all
elliptic elements may be identified with the (topological) fundamental
group
$\pi ^{top}_1(\Gamma )$ of the graph
$\Gamma
$ (this is a general property of graphs of groups). All labelled graphs  
representing $G$ thus have the same {\it first Betti number\/}, denoted
by $b$. 

All homomorphisms from $G$ to a free group with
non-abelian image factor through the quotient map $\theta  :G\to
G/G_{ell}\simeq F_b$ (because in a  free group the
commensurator of any non-trivial element is cyclic). Since $G$ always
maps onto $\Z$, the maximum rank of a free
quotient of $G$ is $\max(b,1)$.

\subhead   The modular homomorphism  $\Delta $ \endsubhead

Let $G$ be a non-elementary GBS group.  The set $\E$ consisting of all
non-trivial elliptic elements is stable under conjugation, elements of
$\E$ have infinite order, and any two elements of $\E$ are
commensurable. These properties yield a homomorphism
$\Delta $ from
$G$ to the multiplicative group of non-zero rationals $\Q^*$, defined as
follows. 

Given $g\in G$, choose any   $a\in\E$. There is a relation
$ga^pg\mi=a^q$, with
$p,q$ non-zero, and we define $\Delta (g)=\frac p q$. As pointed out in [\Krn],
it is easily checked that this is independent of the   choices made ($a$ and the
relation), and defines a homomorphism.  We
call $\Delta (g)$ the {\it modulus\/} of $g$.   Note that $\Delta
\circ\alpha =\Delta $ if $\alpha $ is any automorphism of $G$, because
$\E$ is $\alpha $-invariant. 

Let $H$ be a finite index subgroup of $G$. Any GBS $G$-tree is also a GBS
$H$-tree, so $H$ is a GBS group. The modular homomorphism of $H$ is
the restriction of that of $G$. 

Every elliptic element has modulus $1$, so $\Delta $
factors through the free group $G/G_{ell}\simeq F_b$. In particular, $\Delta $
is trivial when $\Gamma $ is a tree. If
$\Gamma
$ is a labelled graph representing $G$, one has $G/G_{ell}\simeq\pi
^{top}_1(\Gamma )$ and the modulus   may be computed as follows
(see [\BK, \For]): if
$\gamma \in \pi ^{top}_1(\Gamma )$ is represented by an edge-loop
$(e_1,\dots,e_m)$, its modulus is simply
$\ds\prod_{j=1}^m  \frac{\lambda _{e_j}}{\lambda _{\ov e_j}}$.

We denote by $M$ the image of $\Delta $. It is a  
subgroup of $(\Q^*, \times)$. If $G=BS(m,n)$, then $M$ is generated by $\frac
mn$. If $G$ is represented by the labelled graph of Figure 6, then $M$ is
generated by $-1$ and $-\frac32$.

\nom\modu
\thm{Lemma \sta} Let $r=\frac pq$ be  a nonzero rational number,
written in lowest terms. Assume $r\neq\pm1$. 
\roster
\item $r\in M$ if and only if the equation $xy^px\mi=y^q$ has a solution
with
$y\neq1$.
\item If $r\in\Z$, then $r\in M$ if and only if $G$ contains a subgroup
isomorphic to
$BS(1,r)$. 
\endroster
\fthm

\demo{Proof} If $\frac pq \in M$,   the equation
$xy^{np}x\mi=y^{nq}$  has a non-trivial solution for some $n\in\Z$, so
$xy^px\mi=y^q$ has a non-trivial solution. Conversely, if
$|p|\neq |q|$ and 
$xy^px\mi=y^q$  has a non-trivial solution, then $y$ must be elliptic  
and therefore
$\frac pq=\Delta (x)\in M$. We have proved $(1)$.

If $xy^rx\mi=y $ with $y\neq1$ and $r$ an integer different from $-1,0,1$,  
then
 $H=\langle x,y\rangle$ is a solvable GBS group, and $\Delta (x)=r$, so
$H\simeq BS(1,r)$ (one may also show $H\simeq BS(1,r)$  by arguing that    
the only torsion-free proper quotient of $BS(1,r)$ is $\Z$).
\cqfd\enddemo

\example {Remarks} $\bullet$ The values $ \pm1$ are special. If $G$  is
represented by a labelled tree
  $\Gamma $ containing a $(2,2)$-edge, then   $xy x\mi=y\mi$ has a
nontrivial solution because $G$ contains a Klein bottle group, but $-1\notin
M$ because $\Gamma $ is a tree. Conversely,
$1$ always belong to $M$, but
$BS(1,n)$ does not   contain $\Z^2=BS(1,1)$ for $|n|>1$.  

$\bullet$  It is probably not true that $G$
always contains
$BS(p,q)$  if
$\frac pq\neq\pm1$ is a modulus.

$\bullet$ Using $\Delta $, it is easy to show that the isomorphism type of
$BS(m,n)$ determines $m$ and $n$ (normalized by $1\le m\le|n|$) [\Mold]. In
most cases,
$m$ and
$n$ are determined by $m/n$ (given by $\Delta $) and $|m-n|$ (given by
abelianizing). To obtain $m$ from $BS(m,m)$, observe that the quotient of
$BS(m,m)$ by its center is $\Z*\Z/m\Z$. 
\endexample

 We say that $G$ has {\it trivial modulus\/} if $M=\{1\}$ (we often write this as
$\Delta =1$). It is {\it unimodular\/} if
$M\inc\{1,-1\}$, equivalently  if   $xy^px\mi=y^q$
with $y\neq 1$ implies $p=\pm q$. As in [\For], we say that $G$ has {\it no
nontrivial integral modulus\/} if
$M\cap\Z\inc\{1,-1\}$.  This is equivalent to saying that $G$ contains no
solvable Baumslag-Solitar group $BS(1,n)$ with $n\ge2$ (we may take $n>0$
because $BS(1,-n)$ contains $BS(1,n^2)$).

\subhead Unimodular groups \endsubhead

\nom\cent
\thm{Proposition \sta} Let $G$ be a non-elementary GBS group. The
center $Z(G)$ of
$G$ is infinite cyclic if $G$ has trivial modulus, trivial otherwise. It acts as the
identity on any GBS tree.
\fthm

\demo{Proof} Let $T$ be any   GBS  tree (recall that $T$ is always assumed to
be minimal).    If
$a$ is central (more generally, if $\langle a \rangle$ is normal), it is elliptic,
as otherwise its axis would be a
$G$-invariant line. The fixed point set of $a$ is a $G$-invariant subtree,
so equals $T$ by minimality. This shows that
$Z(G)$ is contained in the kernel of the action (elements acting on $T$
as the identity). In particular, it  is trivial or cyclic. 

If $\Delta (g)\neq1$ and $a\in\E$, there is a relation $ga^pg\mi=a^q$
with $p\neq q$, so $a$ cannot be central. This shows that  $Z(G)$ is trivial if $G$
does not have trivial modulus. If $\Delta $ is trivial, choose any finite
generating system
$s_i$ for
$G$, and   $a\in\E$. For each $i$, there is a relation
$s_ia^{n_i}s_i\mi=a^{n_i}$ with $n_i\neq0$. It follows that some power
of $a $ is central. 
\cqfd\enddemo

\thm{Proposition \sta} Let $G$ be a non-elementary   GBS group. The
following are equivalent:
\roster
\item
$G$ is unimodular.
\item
$G$ contains a normal infinite cyclic subgroup $Z$. 
\item
$G$ has a finite index subgroup isomorphic to $F_n\times \Z$ for some
$n>1$.  
\endroster The quotient of $G$ by any normal infinite cyclic subgroup
$Z$   is virtually free. 
\fthm

\demo{Proof} Suppose $G$ is unimodular. The kernel of $\Delta $ has
index 1 or 2, and has trivial modulus. Its center   is
infinite cyclic and characteristic in
$G$, so (1) implies (2). 

Suppose $Z$ is infinite cyclic and normal. Let $T$ be any   GBS
tree. As in the   proof of Proposition \cent, one shows that $Z$ is
contained in the kernel of the action. The quotient
$G/Z$ acts on
$T$ with finite stabilizers, so is virtually free. This easily implies that
$G$ is virtually $F_n\times\Z$. 

If $G$ is virtually $F_n\times\Z$, its modulus is trivial on a finite index
subgroup, so $M$ is finite, hence contained in $ \{1,-1\}$.
\cqfd\enddemo

\example{Remarks} $\bullet$ Let $T$ be a   GBS tree. If $G$ is
unimodular,   we have seen that $ Z(\ker(\Delta ))$ is contained in the
kernel of the action on $T$. Conversely, if the action has a nontrivial  kernel
$K$, then $G$ is unimodular (because $K$ is normal and cyclic), and furthermore
$K=Z(\ker(\Delta ))$. To see this, simply note that, if
$a$ generates $K$, one has $gag\mi=a^{\pm1}$ for any $g\in G$,
so $a$ commutes with $\ker(\Delta )$.

$\bullet$ Two non-solvable  GBS groups are quasi-isometric if and only
if they are both 
 unimodular or both non-unimodular [\Wh].  Any torsion-free group
quasi-isometric to
$F_n\times
\Z$ with $n>1$ is a unimodular GBS group [\MSW].

$\bullet$ A GBS group is residually finite if and only if it is solvable or
unimodular [\Lrang].
\endexample

\subhead   Twisted conjugacy classes \endsubhead 

Let $\alpha :G\to G$ be an endomorphism. Two elements $g,g'\in G$ are
{\it $\alpha
$-conjugate\/} if there exists $h$ such that $g'=hg\alpha (h)\mi$. The
number of
$\alpha $-conjugacy classes is the {\it Reidemeister number\/} of
$\alpha $, denoted by
$R(\alpha )$. It is relevant for fixed point theory  (see [\Fed]).

\thm{Proposition \sta} Let $\alpha :G\to G$ be an endomorphism of a
non-elementary GBS group. If one of the following conditions holds,
then
$R(\alpha )$ is infinite:
\roster \item $\alpha $ is surjective.
\item  $\alpha $ is injective and $G$ is not unimodular.
\item $G=BS(m,n)$ with $|m|\neq|n|$, and the image of $\alpha $ is not
cyclic. 
\endroster
\fthm

This generalizes results of [\Fed] about Baumslag-Solitar groups.

\demo{Proof}   First suppose that $G$ is unimodular and $\alpha $ is
surjective. The group $G$ is residually finite (because it is virtually
$F_n\times\Z$), hence Hopfian. We therefore assume that $\alpha $ is an
automorphism. The subgroup
$Z=Z(\ker(\Delta ))$ is characteristic, so $\alpha $ induces an
automorphism
$\beta $ on the virtually free group $G/Z$. As $G/Z$ is a
non-elementary (word) hyperbolic group,
$R(\beta )$ is infinite [\LL, \Fe]. This implies that $R(\alpha )$ is infinite.

From now on, we assume that $G$ is not unimodular. If $\alpha $ is an
automorphism, we know that $\Delta \circ\alpha =\Delta $, so $\alpha
$-conjugate elements  of $G$ have the same   modulus.  As $M$ is infinite, we
get $R(\alpha )$ infinite. This argument works in the general case, but we
have to prove $\Delta \circ\alpha =\Delta $ for endomorphisms satisfying (1),
(2), or (3).

We first claim that 
$\alpha $ does not factor through $\theta  :G\to G/G_{ell}\simeq F_b$. This
is clear if (2) or (3) holds. If a surjective $\alpha $ factors as $\rho
\circ\theta  $, then  $\theta 
\circ\rho  $   is a non-injective epimorphism from $F_b$ to itself, a
contradiction because free groups are Hopfian.

We can now show   $\Delta \circ\alpha =\Delta $. Since $\alpha $ does
not factor through
$\theta  $,  there is an elliptic 
$a$ with
$\alpha (a)\neq 1$. As $G$ is not unimodular, there is a relation
$g_0a^mg_0\mi=a^n$ with $|m|\neq |n|$. From $\alpha (g_0)\alpha
(a)^m\alpha (g_0)\mi=\alpha (a)^n$, we deduce that $\alpha (a)$ is
elliptic. Thus $\alpha (a)$ is a nontrivial elliptic element, and may be
used to compute $\Delta $. Given any $g$, we have a relation  $g a^pg
\mi=a^q$.  We then write  $\alpha (g )\alpha (a)^p\alpha (g )\mi=\alpha
(a)^q$, showing that $g$ and
$\alpha (g)$ have the same modulus $\frac pq$.
\cqfd\enddemo

\head \sect The automorphism group of a GBS tree  \endhead

\subhead  General facts\endsubhead

Let $G$ be any finitely generated group. As above, we consider
$G$-trees up to equivariant isomorphism. There is a natural action of
$\Out(G)$ on the set of $G$-trees, given by  precomposing an action of
$G$ on $T$ with an automorphism of
$G$ (composing with an inner automorphism does not change the  tree).  

Given $T$, we denote by $\Out^T(G)\inc\Out(G)$ its stabilizer: $\Phi $ is
in
$\Out^T(G)$  if and only if $T$, with the action of $G$ twisted by $\Phi $,
is equivariantly isomorphic to $T$ with the original action.  When $T$ is
irreducible, this is equivalent to saying that the length function $\ell$
satisfies $\ell\circ\Phi =\ell$.

We recall results of [\LGD] about $\Out^T(G)$. We assume that $T$ is
minimal and is not a line  (but there is no condition
on edge and vertex stabilizers in this subsection). The quotient graph of
groups is denoted by
$\Gamma
$, its vertex set by
$V$. 

There is a natural homomorphism from $\Out^T(G)$ to the
symmetry group   of $\Gamma $ (viewed as a graph
with no additional structure). The kernel is a finite index subgroup
$\Out_0^T(G)$, and there is  a homomorphism $\rho
:\Out_0^T(G)\to\prod_{v\in V}\Out(G_v)$. All automorphisms of $G_v$
which occur in the image of
$\rho
$ preserve the set of conjugacy classes of incident edge groups. 

The kernel of
$\rho $ is generated by the group of twists $\T(T)$ together with
automorphisms called bitwists (bitwists belong to $\T(T)$  when vertex
groups are abelian). The group $\T(T)$, which we also denote by $\T(\Gamma
)$, will play an important role in the sequel. Before defining
it, we mention that, when edge groups are cyclic,  there is a further
finite index subgroup 
$\Out_1^T(G)\inc\Out_0^T(G)$ with $\Out_1^T(G)\cap\ker \rho =\T(T)$. It will
be used in Section 7.

To define $\T(T)$, we first consider    an oriented edge  $e$ of $ \Gamma $,
with origin
$o(e)=v$. Let
$G_e$,
$G_v$ be the corresponding edge and vertex groups, with $G_e$
identified to its image in $ G_v$.   We  denote by $Z_{G_v}(G_e)$ the
{\it centralizer\/} of $G_e$ in $G_v$. 

Given $z\in Z_{G_v}(G_e)$,
we define the
{\it twist\/}
$D(z)\in\Out(G)$ by $z$ around $e$ as follows (see [\LGD] for details). If $e$ is
separating, it expresses $G$ as an amalgam
$G=G_1*_{G_e}G_2$.  Then $D(z)$ is defined as the identity on $G_1$, 
and conjugation by $z$ on
$G_2$.    If $e$
does not separate, $G$ is an HNN-extension and $D(z)$ maps the
stable letter $t$ to $zt$ (keeping the base group fixed). 

The {\it group of twists $\T(T)$\/}, or $\T(\Gamma )$,  is   the
subgroup of
$\Out(G )$ generated by all twists. As twists around distinct edges commute,
$\T(\Gamma )$ is a quotient of  
$\prod Z_{G_{o(e)}}(G_e)$, the product being taken over all oriented
edges of $\Gamma $. 
Proposition 3.1 of   [\LGD] says that only two types of relations are
needed to obtain a presentation of $\T(\Gamma )$. 

For each pair   of
opposite edges $(e,\ov e)$, there are   {\it edge relations\/} associated to
elements $z$ in the center $Z(G_e)$ (twisting by $z$   near
the origin of $e$ defines the
same outer automorphism as twisting   by $z\mi$ near the origin  of
$\ov e$). For each vertex $v$, there are {\it vertex relations\/} associated to
elements
$z\in Z(G_v)$ (twisting by $z$ simultaneously around all edges with origin
$v$ defines an inner automorphism). 

\nom\te
\example {Remark \sta}  Let $e\inc\Gamma $ be a segment such that both
adjacent vertex groups are abelian. Using the vertex relations, one sees that
$\T(\Gamma )$ is generated by the groups $Z_{G_{o(f)}}(G_f)$ with $f\neq
e,\bar e$. Collapsing $e$ yields a new graph of groups whose group of twists
contains $\T(\Gamma )$. 
\endexample

Our main tool for finding free groups $F_2$  in $\Out(G)$ will be:

\nom\tw
\thm{Lemma \sta} Let $\Gamma $ be a  minimal graph of groups,
with fundamental group $G$. Let $e$ be an edge
with origin
$v$, and  let $G_e$, $G_v$ be the corresponding groups. The subgroup
$\T(\Gamma )\inc \Out(G)$  maps onto
$Z_{G_v}(G_e)/\langle Z( G_v),Z(G_e)\rangle$.  
\fthm

We  denote by   $\langle Z(
G_v),Z(G_e)\rangle$ the (obviously normal) subgroup
generated by the centers of $G_v$ and $G_e$.

\demo{Proof} 
 Divide $\T(\Gamma )$ by (the image of) all factors $Z_{G_{o(f)}}(G_f)$ for
$f\neq e$ (including $f=\ov e$). The only relations which remain are those
involving
$Z_{G_v}(G_e)$, namely  edge relations associated to $(e,\ov e)$ and
vertex relations associated to $v$. The quotient is precisely 
$Z_{G_v}(G_e)/\langle Z( G_v),Z(G_e)\rangle$.  
\cqfd\enddemo

\subhead  The group of twists  of a GBS tree\endsubhead

Now let $G$ be a non-elementary GBS group. We consider the  action of
$\Out(G)$ on the set of GBS trees.  The corresponding action   on the
set of marked graphs is by changing the marking.  If $T$ is
rigid, then
$\Out^T(G)=\Out(G)$.

The  group of twists $\T(T)$ is a finitely generated   abelian group.
The presentation recalled above may be rephrased as follows (we  
use additive notation). 

Given an oriented edge $e$ of $ \Gamma $, there is one generator
$D_e$.  If $e$ is separating,  $D_e$ is   the identity on $G_1$, and
conjugation by $x_v$ on
$G_2$ (with $x_v$ the  generator of the vertex group
at $v=o(e)$, and   $G=G_1*_{G_e}G_2$ as above).  If $e$
does not separate, choose a maximal tree $\Gamma _0$ not containing
$e$. In the corresponding presentation of $G$, define $D_e$ as mapping
$t_e$ to $x_{o(e)}t_e$ and keeping all other generators fixed. 

In terms of these generators $D_e$, the relations are the following. 
For each pair  of opposite
edges $(e,\ov e)$, there is an edge relation
$\lambda _e D_e+\lambda _{\ov e} D_{\ov e}=0$, implied by the relation
$(x_{o(e)})^{{\lambda (e)}}=  (x_{o(\ov e)})^{{\lambda (\ov e)}}$ or 
$t_e(x_{o(e)})^{{\lambda (e)}}t_e\mi=   (x_{o(\ov e)})^{{\lambda (\ov e)}}$. For
each vertex $v$, there is a vertex relation $\sum_{e\in E_v}D_e=0$, with
$E_v$ the set of edges with origin $v$.  

\example{Remark} The group $G$ is in a natural way the fundamental group of
a 2-complex consisting of annuli (corresponding to edges of $\Gamma $) glued
to circles (corresponding to vertices). One can consider the subgroup
$DT(\Gamma )$ of
$\T(\Gamma )$ generated by Dehn twists supported in the annuli. It is easy
to see that it has finite index. One may also show
$\T(\Gamma )=DT(\Gamma ')$, where $\Gamma '$ is  a (non reduced)
graph obtained from
$\Gamma
$ by elementary expansions.    
\endexample

 Recall that 
$b $ is the first Betti number of any labelled graph $\Gamma $
representing
$G$. 

\nom \kk
\thm{Proposition \sta} Let $G$ be a non-elementary GBS group. Define $k$ as
$b$ if $G$ has trivial modulus,
$b-1$ if not. 
\roster \item The torsion-free rank of the abelianization $G_{ab}$ is
$k+1$.
\item Let $\Gamma $ be any labelled graph representing $G$. The
torsion-free rank of the abelian group $\T(\Gamma )$ is
$k$.  

\item If $\Gamma _0\inc\Gamma $ is a 
maximal subtree, the twists $D_e$   around the
 edges of $\Gamma \setminus\Gamma _0$ generate a finite index subgroup of
$\T(\Gamma )$.
\endroster
\fthm

\demo{Proof} Killing all elliptic elements produces an epimorphism $\theta
:G\to F_b$ (see Section 2), so   
$rk(G_{ab})\ge b$. If $\Delta $ is nontrivial, any elliptic element $a$ satisfies a
relation  $ga^pg\mi=a^q$ with
$p\neq q$, so is  mapped trivially to
torsion-free abelian groups. This shows $rk(G_{ab})= b$ in this case. 

If $\Delta $ is trivial, fix $\Gamma $ and $\Gamma _0$. It
follows from the presentation of $G$ given earlier that $G_{ab}$ is the direct
sum of
$\Z^b$ with the abelian group
$G'$ defined by the following presentation: there is one generator $x_v$
for each vertex of $\Gamma $, and one relation $\lambda _e x_{o(e)}=\lambda
_{\ov e} x_{o(\ov e)}$ for each pair of opposite edges
$(e,\ov e)$. We show that $G'$  maps non-trivially to $\Z$. It is easy to map the
generators
$x_v$ to
$\Z$ in such a way that relations associated to edges in $\Gamma _0$ are
satisfied. Using the formula
$ \Delta (\gamma )=\prod_{j=1}^m  \frac{\lambda _{e_j}}{\lambda _{\ov
e_j}}$ (see Section 2), one sees that  the remaining relations are  
automatically satisfied. We get $rk(G_{ab})= b+1$.

Assertion (2) follows immediately from Proposition \cent{} and the exact
sequence 
$$0\to Z(G)\to \Z^{\nu +\zeta }\to\Z^{2\zeta }\to\T\to0$$
given by Proposition 3.1 of  [\Ler], where $\nu $ (resp\. $\zeta $) is the number
of vertices (resp\. edges) of $\Gamma $ (we are grateful to M\. Clay for
suggesting this short argument).

Assertion (3) follows from the presentation of $\T(\Gamma )$ in terms of the
generators $D_e$. If we add the relations $D_e=D_{\ov e}=0$ for
$e\notin\Gamma _0$, the quotient is the group of twists associated to the
labelled graph $\Gamma _0$, so is finite by Assertion (2). 
 \cqfd\enddemo

\nom\foo
\example{Remark  \sta}
One may decide whether a given
$D_e$ has finite or infinite order. If $e$ does not separate, the order of
$D_e$ is finite if and only if $G$ has nontrivial modulus, but every loop
not containing $e$ has trivial modulus. If $e$   separates, the order
is infinite if and only if each component of $\Gamma \setminus\{e\}$
contains a loop with nontrivial modulus. 
\endexample

 The   groups
$\T(T)$ associated to different GBS trees are abstractly commensurable by
Proposition
\kk. We show that they are commensurable as subgroups of $\Out(G)$.

\nom\commens
\thm{Proposition \sta} If $T,T'$ are two GBS  trees, then  $\T(T)$ and 
$\T(T')$ are commensurable subgroups of $\Out(G)$. 
\fthm

\demo{Proof} By Corollary \defor, it suffices to show that $\T(T')$ is
commensurable with $\T(T)$  if 
$T'$ is obtained from
$T$ by an elementary collapse.  Consider the corresponding graphs
$\Gamma ,\Gamma '$. Let $e=vw\inc \Gamma $ be the collapsed edge.
We  assume $\lambda _e=1$, and we denote
$\lambda_{\ov e}$ by $\lambda $, so $G_v$ has index $\lambda $ in
$G_w$.  

 Let $F$ be the set of oriented edges of $\Gamma $ other than $e,\ov e$.
The group
$\T(T)$ is the subgroup of $\Out(G)$ generated by the twists $D_f$,
$f\in F$ (see Remark \te). Similarly, $\T(T')=\langle
D'_f\mid f\in F\rangle$, as $F$ may be viewed as  the set of
oriented edges of
$\Gamma '$. Moreover, we have $D_f=D'_f$ if the origin of $f$ is not $v$,
and
$D_f=\lambda D'_f$ if it is because the collapse replaces the vertex
group $G_v$ of
$\Gamma $ by the larger group $G_w$. This shows that $\T(T')$ contains
$\T(T)$ as a subgroup of finite index. 
\cqfd\enddemo

\nom\divi
\example{Remark \sta} The index of $\T(T)$ in $\T(T')$
divides a power of the label $\lambda $. This will be used in the proof of
Theorem 3.12.
\endexample

\subhead  Applications  \endsubhead

We apply the preceding results to the study of $\Out^T(G)$, using the
following fact:

\thm{Proposition \sta{}  [\LGD]}  $\T(T)$ has finite index in  $\Out^T(G)$.
\fthm

This follows from Theorem 1.6 of  [\LGD], as edge and vertex groups
have finite outer automorphism groups. More precisely, let us show:

\nom\unib
\thm{Proposition \sta} Given $G$, the index of $\T(T)$ in $\Out^T(G)$ is
uniformly bounded (independently of
$T$). 
\fthm

\demo{Proof}
Consider the chain of
subgroups
$\T(T)\inc
\ker\rho
\inc\Out_0^T(G)\inc\Out^T(G)$ mentioned at the beginning of this
section. We check
that each group has uniformly bounded index in the next. 

The index of
$\Out_0^T(G)$ in
$\Out^T(G)$ is   bounded by the order of the symmetry group of $\Gamma $.
The number of edges of $\Gamma $ is not always uniformly bounded, but the
first Betti number  
  is fixed, and there is a uniform bound
  for the number $d$ of terminal vertices  
(because  adding the relations $x_v=1$, for
$v$ non-terminal, maps $G$ onto the free product of
$d$ nontrivial finite cyclic groups). This is enough to bound the symmetry
group. 

 The map $\rho $ describes how automorphisms act on vertex groups.
Since these groups are all  
commensurable, and isomorphic to $\Z$, the image of
$\rho $ has order at most 2, so $\ker\rho $ has index at most 2 in
$\Out_0^T(G)$. Finally,
$\ker\rho $ is generated by $\T(T)$ together with bitwists. As vertex groups
are abelian, bitwists belong to $\T(T)$, so
$\T(T)=\ker\rho $.
\cqfd\enddemo

If $T $ is rigid, then $\Out^T(G)=\Out(G)$. We get:

\thm{Theorem \sta} If   $G$ is algebraically rigid, then $\Out(G)$
contains  
$\Z^k$ as a subgroup of finite index.\cqfd
\fthm

In general, we have:

\nom\comm
\thm{Theorem \sta} Up to commensurability, the subgroup $\Out^T(G)$
of $\Out(G)$ does not depend on
$T$.  It contains
$\Z^k$ with finite index.   \cqfd
\fthm

Another proof of the first assertion (and therefore of Proposition
\commens, but not of Remark \divi) is given in Section 5. Also note the
following related result: 

\nom\clay
\thm{Theorem \sta{} [\Clc]} Any subgroup of $\Out(G)$
commensurable with  a subgroup of
$\Out^T(G)$ is contained in $\Out^{T'}(G)$ for some GBS tree $T'$. \cqfd
\fthm

We now prove:

\nom\tooo
\thm{Theorem \sta} The set of prime numbers $p$ such that $\Out(G)$
contains non-trivial $p$-torsion is finite. 
\fthm

\demo{Proof} As any torsion element of $\Out(G)$ is contained in some
$\Out^T(G)$ by [\Clc], and the index of $\T(T)$ in $\Out^T(G)$ is
uniformly bounded by Proposition \unib, it suffices to control torsion in  
groups of twists. 

First note that the set of prime numbers dividing a label of $\Gamma $
does not depend on $\Gamma $, as it does not change during an
elementary collapse.   Call it
$\P$. If
$T$ and
$T'$ are related by an elementary collapse, Remark \divi{} shows that 
$\T(T)$ and
$\T(T')$ have torsion at the same primes, except possibly those in $\P$.
This implies that only finitely many primes may appear in the torsion
of a group of twists: those in
$\P$, and those in the torsion of $\T(T_0)$ for some fixed $T_0$. 
\cqfd\enddemo

We have seen that $\Out(BS(2,4))$ contains arbitrarily large $2$-torsion. The
proof of Theorem \tooo{} also shows:

\thm{Corollary \sta} If $p$ is a prime number such that $\Out(G)$
contains $p$-torsion of arbitrarily large order, then $p$ divides at
least one label of each labelled graph representing $G$. \cqfd
\fthm

We do not know whether $p$ must divide some integral modulus. 

\head  {\sect Unimodular groups}\endhead

Let $\Gamma $ be a labelled graph representing a non-elementary
unimodular group $G$, and $T$ the associated Bass-Serre tree. We denote by
$G^+$   the kernel of
$\Delta:G\to\{\pm1\}
$ ({\it positive\/} elements). All elliptic elements are positive. Let $Z$ be the
center of
$G^+$. We know that it is cyclic, characteristic in $G$, and acts as the identity
on 
$T$.   

We fix a nontrivial $\delta \in Z$. If $\Delta =1$ we may take $\delta  $ to be a
generator $\delta _0$, but the study of $\Out(G)$ when $\Delta \neq 1$ will
require 
$\delta $ to be $(\delta _0)^4$. Note that any
generator of an edge or vertex group is a root of $\delta
$.

 Let $Z'$ be the cyclic group
generated by $\delta $.
 There is an exact sequence
$\{1\}\to Z'\to G\to H\to\{1\}$ with  $H$
virtually free.
The group $H$ is the fundamental group of a   graph of groups with  the same
underlying graph. Vertex and edge groups are finite cyclic groups,
the  order being the index of
$\langle\delta \rangle$ in the original group. We denote by $\bar g$ the image
of $g\in G$ in $H$.

Since $Z'$ is characteristic in $G$, there are   natural homomorphisms
$ \Aut(G)\to \Aut(H)$ and $ \Out(G)\to \Out(H)$. The basic example is
$\Out(F_n\times\Z)$, which contains the semi-direct product
$\Z^n\rtimes\Out(F_n)$ with index $2$ (the factor $\Z^n$ should be thought of
as $\text{Hom}(F_n,\Z)$). But the following examples illustrate a few of the
subtleties involved when trying to   lift automorphisms from
$H$ to  
$G$.

\example{Examples} 

$\bullet$ Let $G$ be $\langle a,b\mid a^3=b^3\rangle $ and $H$ be $\langle \bar
a,\bar b\mid
\bar a^3=\bar b^3=1\rangle $. The automorphism of
$H$ mapping $\bar a$ to $\bar a\mi$ and $\bar b$ to $\bar b$ does not
lift to $G$.

$\bullet$ $G$ is $BS(3,3)=\langle a,t\mid ta^3t\mi=a^3\rangle $ and $H$ is
$\langle
\bar a,\bar t\mid
\bar  a^3=1\rangle $. The automorphism   fixing $\bar a$ and sending $\bar t$
to
$\bar t\bar a$   has order 3,   but all its lifts have infinite order.

$\bullet$ $G$ is $BS(2,-2)=\langle a,t\mid ta^2t\mi=a^{-2}\rangle $ and $H$ is
$\langle
\bar a,\bar t\mid \bar a^2=1\rangle $. Conjugation by
$\bar a$ in
$H$ has  lifts of order 2, such as $a\mapsto a, t\mapsto  ata$, or
$a\mapsto a\mi,t\mapsto a ta\mi$, but no lift of order 2 is inner.

\endexample

Let $H^+$ be the image of
$G^+$ in
$H$. If $\Delta \neq1$,  it has index 2 (because $\delta $ is positive).
There are only finitely many conjugacy classes of torsion elements in
$H$ (they all come from vertex groups). All torsion  elements of $H$
belong to
$H^+$, but a conjugacy class in $H$ may split into two classes in $H^+$. 

We shall now define  a homomorphism $\tau  :G \to
\text{Isom}(\R)$ (it is similar to the homomorphism $G'\to\Z$ constructed in
the proof of Proposition \kk). We fix a maximal tree
$\Gamma _0\inc\Gamma $. Recall the presentation of
$G$ with generators $x_v,t_\varepsilon $ and relations of the form
$x_v^m= x_w^n$ or   $t_\varepsilon x_v^mt_\varepsilon \mi= x_w^n$. 

To define $\tau  $, send
$\delta $ to $x\mapsto x+1$, send $x_v$ to $x\mapsto x+1/n_v$ if
$\delta=x_v^{n_v}$, send
$t_\varepsilon $ to
$x\mapsto \Delta (t_\varepsilon )x$, and check that the relations are
satisfied. This
$\tau  $ is not canonical (it depends on the choice of  
$\Gamma _0 
$); it is uniquely defined on elliptic elements once  
$\delta $ has been chosen.

 The  image of $\tau  $ in $\text{Isom}(\R)$ is infinite cyclic if
$\Delta =1$, infinite dihedral if $\Delta \neq1$. Its kernel contains no
nontrivial elliptic element. The coefficient of $x$ in $\tau  (g)$ is $\Delta
(g)$, and
$\tau  (gcg\mi)=\tau  (c)^{\Delta (g)}$ if
$c$ is positive (in particular if $c$ is elliptic).

The map $\tau  $ induces a map $\ov\tau  $ from $H$ to a finite group $F$
(the quotient of the image of $\tau  $ by $x\mapsto x+1$).  The group $F$
is cyclic if
$\Delta =1$, dihedral if $\Delta \neq 1$.

\definition{Definition} We define the finite index 
subgroup $\Aut_0(H)\inc\Aut(H)$ as the set of automorphisms $\ov\alpha $
such that:
\roster
\item
$ \ov\alpha (H^+)=H^+$;
\item
 $\ov\alpha $ acts trivially on the set of
$H^+$-conjugacy classes of torsion elements;
\item
 $\ov\tau 
\circ\ov\alpha =\ov\tau  $.  
\endroster
\enddefinition

\nom\liftt
\thm{Lemma \sta}  Let $\ov\alpha \in\Aut_0(H)$. There exists a
unique  lift 
 $\alpha \in\Aut(G)$ such that  
$\tau  \circ\alpha  =\tau  $. It satisfies $\alpha (\delta )=\delta $.
\fthm

\demo{Proof} 
 Uniqueness is easy: $\alpha (g)$ is determined up to a power of
$\delta $, and that power is determined by applying $\tau  $.

We define $\alpha $ on the generators of $G$. 
  In $H$, the element $\bar x_v$ has finite order and
therefore is mapped by $\ov\alpha $ to  
$\bar g_v   \bar x_v  \bar g_v \mi$ for some $\bar g_v\in H^+$. We define
$\alpha (x_v)=g_vx_vg_v\mi$, where
$g_v\in G^+$ is any lift of $\bar g_v$. Note that $\tau  (\alpha (x_v))=\tau 
(x_v)$ because $x_v$ and $g_v$ are positive. If  $x_v^{m} =x_w^{n}  $
is a relation, then $\alpha (x_v)^{m} \alpha (x_w)^{-n }$ is $1$ because
it is killed both in
$H$ and by $\tau 
$. Note that $\alpha (\delta )=\alpha (x_v^{n_v})=g_vx_v^{n_v} g_v\mi=\delta
$.

Now consider a generator $t_\varepsilon $, and a lift
  $u_\varepsilon $   of $\ov\alpha (\bar t_\varepsilon )$.
 Since  $\ov\tau 
\circ\ov\alpha =\ov\tau  $, the elements 
$t_\varepsilon $ and $u_\varepsilon $ have the same image in $F$, so $\tau 
(t_\varepsilon u_\varepsilon  \mi)$ is   translation by an integer
$n_\varepsilon $. We define $\alpha (t_\varepsilon )$ as
$ \delta {}^{n_\varepsilon } u_\varepsilon $, so that $\tau  (\alpha
(t_\varepsilon ))=\tau  (t_\varepsilon )$.  Given a   relation  
$t_\varepsilon x_v^{m}t_\varepsilon \mi=x_w^{n }$, the relation
$\alpha (t_\varepsilon )\alpha (x_v)^{m}\alpha (t_\varepsilon )\mi=\alpha
(x_w)^{n }$ holds  modulo $\delta $. It also holds when we apply $\tau  $, so it
holds in $G$.

 We have constructed an endomorphism of $G$ fixing $\delta $ and
inducing $\ov\alpha $, and this forces it to be an automorphism.  \cqfd

\enddemo

Let $\Aut_0(G)\inc\Aut(G)$ be the finite index subgroup consisting of
automorphisms fixing
$\delta
$ and mapping into
$\Aut_0(H)$. We know that the map $\varphi :\Aut_0(G)\to\Aut_0(H)$
is onto and has a section. We consider its kernel.

\nom\kernn
\thm{Lemma \sta} The kernel $N$ of
$\varphi  :\Aut_0(G)\to\Aut_0(H)$ is isomorphic to
$\Z^{b  }$.  It is generated by twists by $\delta $ around the edges of
$\Gamma
\setminus \Gamma _0$.
\fthm

Recall that $b$ is the first Betti number of $\Gamma $.

\example{Remark} It is a general fact that, whenever $Z\inc G$ is
characteristic, the kernel of the map
$ \Aut(G)\to \Aut(Z)\times\Aut(G/Z)$ is abelian [\RV, Proposition 2.5].
To see this, take $\alpha _1,\alpha _2$ in the kernel. Write
$\alpha _1(g)=z_1g$ and
$\alpha _2(g)=gz_2$ (with $z_1,z_2\in Z$, depending on $g$), and
deduce  $\alpha _1\alpha _2(g)=\alpha _2\alpha _1(g)=z_1gz_2$ (one can
also prove that $z_1$ must be in the center of $G$)
\endexample

\demo{Proof}  Suppose $\alpha \in N$. We have $\alpha (\delta )=\delta
$. If $x$ is a root of $\delta $,   we have  $\alpha
(x )=x  
\delta {}^p$ and
$\delta =x ^{q}$, so that $\delta =\alpha (x ^{q})=\delta \delta
{}^{pq}$ and
$p=0$.   Therefore
$\alpha
$ fixes every elliptic element. 
 Furthermore $\alpha (t_\varepsilon )=\delta {}^{n_\varepsilon
}t_\varepsilon $ for some $n_\varepsilon \in\Z$, so $\alpha $ is a product of
powers of twists by $\delta $. Conversely, each choice of integers
$n_\varepsilon $ determines an automorphism fixing all elliptic elements and
belonging to   $N$.
\cqfd
\enddemo

We have proved:

\nom\semidaut
\thm{Theorem \sta} If $G$ is non-elementary and unimodular, there is a split 
exact sequence
$$\{1\}\to
\Z^{b }\to\Aut_0(G)\overset\varphi  \to\rightarrow\Aut_0(H)
\to\{1\},$$  where $H$ is virtually free and 
 $\Aut_0$ has finite index in $\Aut$. \cqfd
\fthm

We shall now show:

\nom\esout
\thm{Theorem \sta} If $G$ is non-elementary and unimodular, there  is a split
exact sequence
 $$\{1\}\to
\Z^{k}\to\Out_0(G)\overset\psi \to\rightarrow\Out_0(H)
\to\{1\},$$ where $H$ is virtually free and 
 $\Out_0$ has finite index in $\Out$.
\fthm

See Proposition 
\kk{} for the definition and properties of $k$. 

Since  $\Out(H)$ is VF [\KV], this implies:

\thm{Corollary \sta} $\Out(G)$ and $\Aut(G)$ are virtually torsion-free and VF
(they have finite index subgroups with finite classifying spaces).
\cqfd
\fthm

\demo{Proof of Theorem \esout} We denote by $\Out_0$ the image of
$\Aut_0$ in
$\Out$ (note that $\Aut_0$ does not contain all inner automorphisms if
$\Delta \neq1$), and by  
$\hat N$ the image of
$N$ in $\Out(G)$.
Let
$\psi :\Out_0(G)\to\Out_0(H)$ be the natural map.  Note that $\hat N$ is
contained in $\ker\psi $, and has torsion-free rank $k$ by Lemma
\kernn{} and Assertion (3) of Proposition \kk. We shall show $\ker\psi =\hat
N\simeq\Z^k$. We write
$i_g$ for conjugation by $g$.

First assume $\Delta =1$. Then $k=b$, and $\hat N\simeq\Z^b$ because it has
torsion-free rank $b$ and is a quotient of $N\simeq\Z^b$. Since the image of
$\tau   $ is abelian, every conjugation $i_g$ in $G$ satisfies $\tau  \circ
i_g=\tau 
$, and Lemma \liftt{} lifts $i_h\in\Aut_0(H)$ to $i_g$, where $i_g$ is any lift of
$h$. Thus $\psi $ has a section. 

There remains to show $\ker\psi \inc\hat N$.  If
$\alpha \in\Aut_0(G)$ represents an element of $\ker\psi$, its image
$\ov\alpha $ in $\Aut_0(H)$ is conjugation by some $h\in H$. Lift $i_h$ to
$i_g\in\Aut_0(G)$, and consider $i_g\mi\alpha $. It belongs to $N$, and  
has the same image as $\alpha $ in
$\Out(G)$. This implies $\ker\psi =\hat N$.

Now   suppose $\Delta \neq1$. In this case we have to choose $\delta =(\delta
_0)^4$, where $\delta _0$ is  a generator of $Z$ (the center of $G^+$).  We first
show that
$\ker\psi =\hat N$. The argument is the same as before, but we have to prove
that
$i_h$ has a lift $i_g\in\Aut_0(G)$ (we will see that Lemma \liftt{} lifts inner
automorphisms  to inner automorphisms, but we cannot claim it at this point).
Since
$\ov\tau 
\circ\ov\alpha =\ov\tau  $, the image $\ov\tau  (h)$ is central in $F$. Our
choice of $\delta $ ensures that the center of $F$ has order $2$ (it is
generated by the image of $x\mapsto x+1/2$). In particular, $h$ is
positive. If $g\in G$ is a lift of $h$, it commutes with $\delta $ and
therefore   $i_g$ belongs to 
$\Aut_0(G)$.

We now prove   $\hat N\simeq
\Z^{b  -1}$. Recall that 
$\hat N$ 
 has torsion-free rank ${b -1}$. 

Consider the twist $D$ by $\delta  =(\delta_0) ^4$ around the edges
$\varepsilon $ of
$\Gamma
\setminus \Gamma _0$ such that $\Delta (t_\varepsilon )=-1$  (it fixes the
generators $x_v$, and maps $t_\varepsilon $ to $t_\varepsilon$ if $\Delta
(t_\varepsilon )=1$, to  
$\delta t_\varepsilon $ if $\Delta (t_\varepsilon )=-1$).
Note that $D$ is   conjugation by
$(\delta _0)^{2}$. Indeed, the  
$t_\varepsilon $'s with modulus $1$, and the
$x_v$'s, are fixed by $D$ and commute with $\delta _0$, whereas
$(\delta_0) ^4t_\varepsilon =(\delta _0)^{2} t_\varepsilon (\delta _0)^{-2}$ if 
 $t_\varepsilon
\delta t_\varepsilon \mi=\delta
\mi$. Since $D$ belongs
to a basis of
$N$ (see  Lemma \kernn), the image $\hat N$ of  $N$ in
$\Out(G)$ is isomorphic to
$\Z^{b  -1}$. 

Finally, we   show that Lemma \liftt{} lifts inner automorphisms  to
inner automorphisms (and therefore $\psi $ has a section).

Suppose   $i_h$ belongs to $\Aut_0(H)$. Then $\ov\tau  (h)$ is central in
$F$, and therefore $\ov\tau  (h^2)$ is trivial. If $g\in G$ is a lift of
$h^2$, then
$\tau  (g)$ is an integral translation and we can redefine $g$ (multiplying
it  by a power of
$\delta $) so that
$\tau  (g)$ is trivial. Then $\tau  (gug\mi)=\tau  (u)$ for every $u\in G$,
showing that  Lemma \liftt{} lifts conjugation by $ {h^2}$   to
conjugation by $g$ in
$G$. Now
  consider the lift $\alpha $ of  $i_h$ given by Lemma \liftt{}. It satisfies
$\alpha ^2=i_g$, and its image
  in 
$\Out(G)$    belongs to   $\ker\psi $. Since
$\ker\psi =\hat N$ is torsion-free, we conclude that $\alpha $ is inner.
\cqfd

\enddemo

\head \sect The deformation space \endhead

Let $G$ be a non-elementary GBS group. In this section, we work with
{\it metric\/} GBS trees: $T$ is a metric tree, and $G$ acts by
isometries.  Metric trees are considered up to $G$-equivariant
isometry.

Let $\D$ be the space of metric GBS  trees, and $P\D$ its
projectivization (obtained by identifying two trees if they differ by
rescaling the metric). We call $P\D$ the (canonical) {\it
projectivized deformation space\/} of
$G$. 

Choosing a $G$-invariant metric on a given simplicial tree amounts to
assigning a positive length  to each edge of
$\Gamma =T/G$.
This makes 
$P\D$ into a complex. An open  simplex is the set of
trees with a given underlying simplicial tree, a closed simplex is the
set of GBS trees that may be obtained from trees in an open simplex
by collapse moves (closed simplices    have ``faces at
infinity'', as the length of a non-collapsible edge is not allowed to be 0). Every
closed simplex contains reduced trees. 

The group $\Out(G) $ acts on $P\D$. There is a bijection between the set
of orbits of open simplices and the set of (unmarked) labelled graphs
representing $G$ (up to admissible sign changes). 
Standard techniques
show that
$P\D$ has a natural $\Out(G)$-equivariant deformation retraction onto
a simplicial complex (see [\CV, \McM]).

GBS trees are locally finite. This implies that the complex $P\D$ is {\it
locally finite\/}.   Indeed, closed simplices containing $T$ consist of
simplicial trees obtained from $T$ by   expansion moves.  Performing
such moves  on
$T$ amounts to blowing up each vertex
$v$ of
$T$ into a subtree. Since $v$ has finite valence, there are only finitely
many ways of expanding (not taking the metric into account). 
As remarked by M\. Clay, this local finiteness gives another
proof of the first assertion of Theorem \comm.

In general, there are several ways to define a topology on spaces of
trees (equivariant Gromov-Hausdorff topology, axes topology, weak
topology), but because of local finiteness  they all coincide on $P\D$
(see the discussion in [\GL, \GLL]). Clay [\Cld] proved that {\it $P\D$ is
contractible\/} (see also [\GLL]).  By Theorem \comm{}, stabilizers for
the
  action of $\Out(G)$ on $P\D$ are virtually $\Z^k$. 

To sum up: 

\nom\comp
\thm{Proposition \sta}   
$\Out(G)$ acts on the locally finite, contractible, complex $P\D$ with
stabilizers virtually $\Z^k$. \cqfd
\fthm

If $G$ is algebraically rigid, the unique reduced GBS tree  belongs to every
closed simplex,   $P\D$ is a finite complex, and the action of $\Out(G)$ on
$P\D$ has a fixed point.   If $G$ is not algebraically rigid,  we
will show that
$\Out^T(G)$ always has infinite index in
$\Out(G)$ (see Theorem 8.5). All $\Out(G)$-orbits are therefore infinite.

Under suitable hypotheses, we now  show that
$P\D$ is ``small'' and we deduce information on $\Out(G)$.

\subhead   Groups with no nontrivial integral modulus \endsubhead

Suppose that $G$ has no integral modulus other than $\pm1$ 
(equivalently, $G$ does not contain a solvable Baumslag-Solitar group
$BS(1,n)$ with $n>1$). In this case, there are only finitely many
$\Out(G)$-orbits of simplices consisting of reduced trees [\For,
Theorem 8.2], and therefore   $\Out(G)$ acts on the complex $P\D$ with
only finitely many orbits. This implies:

\thm{Theorem \sta} Let $G$ be a non-elementary GBS group with no
integral modulus other than $\pm1$. 
\roster
\item $\Out(G)$ is $F_\infty$ (it has a $K(\pi ,1)$ with finitely many cells
in every dimension).  
\item There is a   bound for the cohomological dimension of
torsion-free subgroups of $\Out(G)$.
\item If $\Out(G)$ is virtually torsion-free, it has a finite index
subgroup with a finite classifying space.
\item $\Out(G)$ contains only finitely many conjugacy classes of finite
subgroups.  
\endroster
\fthm

\demo{Proof} The first three assertions follow from Proposition
\comp{} by standard techniques. Unfortunately, we do not
know whether $\Out(G)$ must be virtually torsion-free when $G$ is not
unimodular. Assertion (4) follows from   Theorems \comm{} and  \clay: any
finite subgroup is contained in some $\Out^T(G)$, and it is well-known that
there are finitely many conjugacy classes of finite subgroups in  a group
which is virtually
$\Z^k$ (a proof appears in [\LGD]).
\cqfd\enddemo 

\subhead   Groups with no strict ascending loop   \endsubhead

It is shown in [\GLL] that, if no reduced labelled graph contains a strict
ascending loop,   there is an
$\Out(G)$-equivariant deformation retraction from $P\D$ onto a
finite-dimensional subcomplex. This implies:

\thm{Theorem \sta}
 If no reduced labelled graph representing $G$ contains a strict
ascending loop, there is a   bound for the cohomological
dimension of torsion-free subgroups of $\Out(G)$. \cqfd
\fthm

\head \sect Free subgroups in $\Out(G)$ \endhead

Let $\Gamma $ be a  labelled  graph. Recall that we consider it up
to admissible  sign changes. In
particular, when we focus on an edge, we will always assume that it is a
$(p,q)$-segment with $p,q\ge1$ ($\ge2$ if $\Gamma $ is reduced), or a
$(p,q)$-loop with $1\le p\le|q|$. Recall that a strict ascending loop is a
$(1,q)$-loop with
$|q|\ge2$. A   pseudo-ascending loop is a
$(p,q)$-loop with $p|q$.

\midinsert
\centerline 
{\includegraphics[scale=.5]
{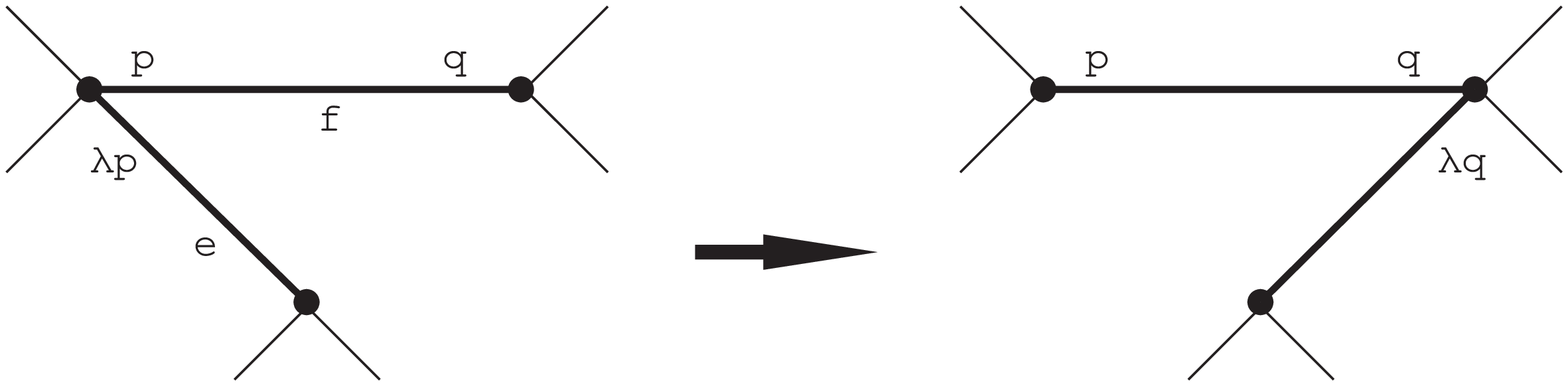}} 
\captionwidth{220pt}
\botcaption
 {Figure \the\fihno}
{Slide move. }
\endcaption
\endinsert

If $e,f$ are distinct oriented
 edges with the same origin $v$, and the label $\lambda _f$ of $f$ at
$v$ divides the label $\lambda _e$ of $e$, one may {\it slide\/} $e$
across
$f$, replacing its label by $\frac{\lambda _e}{\lambda _f}\lambda _{\ov
f}$  (see Figure
\fig); both $e$ and
$f$ may be loops, but they have to be distinct geometric edges ($f\neq
\ov e$). See [\For] for details about slide moves. 
The important thing for
us here is that performing a slide move on  a labelled graph gives another
labelled  graph representing the same group $G$
  (only   the GBS tree changes).

An edge $f$ is a {\it slid\/} edge if some other edge may slide across
$f$  or $\bar f$ (we usually think of slid edges as non-oriented edges).  For
example, any
$(1,q)$-loop is a slid edge if $\Gamma $ contains more than one edge (i.e\. if
$G$ is not solvable).  

\midinsert
\centerline 
{\includegraphics[scale=.45] 
{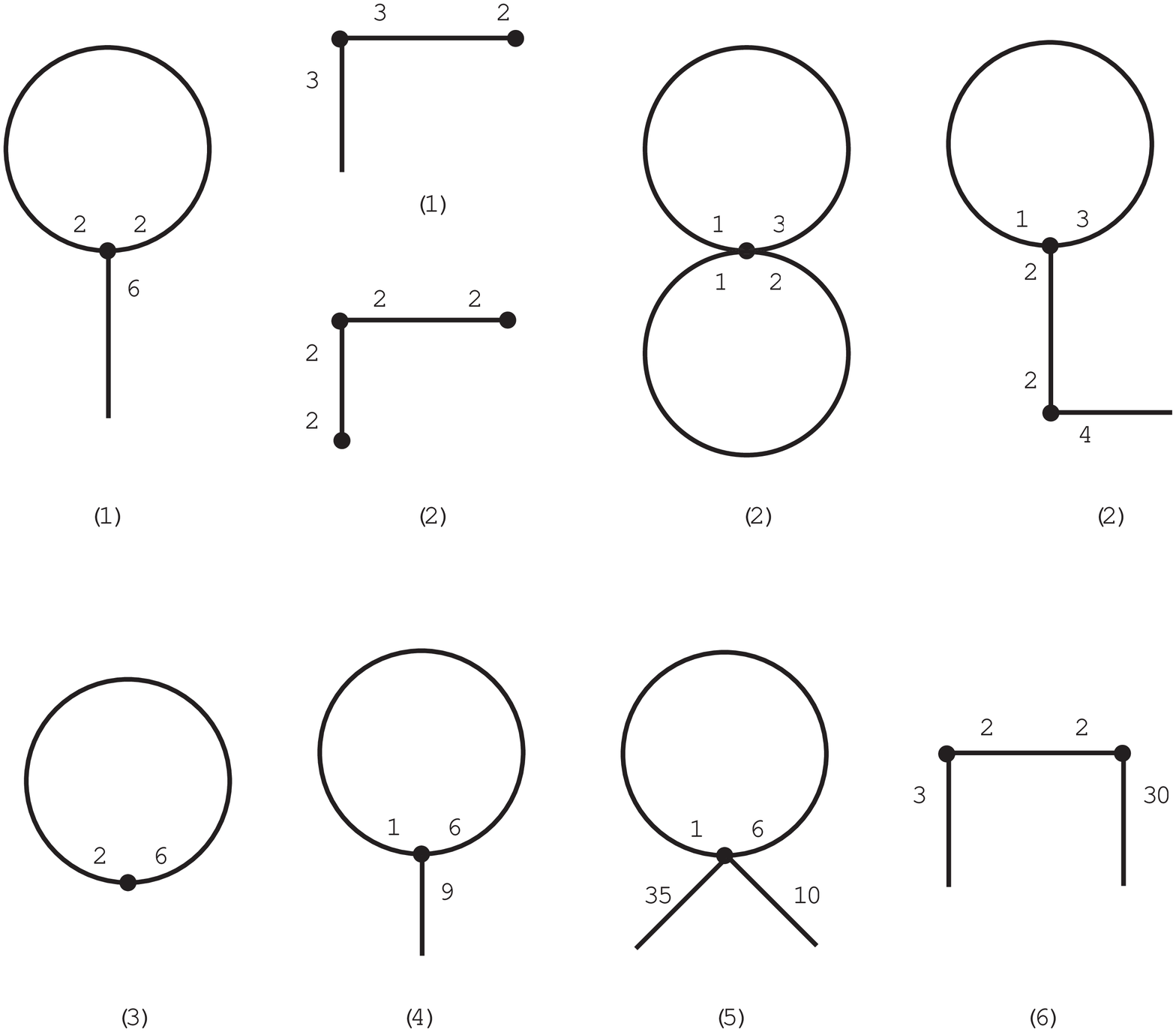}}
\botcaption
 {Figure \the\fihno}
{If   $ \Gamma $  contains one of these graphs,   $\Out(G)$ contains
$F_2$.} 
\endcaption
\endinsert

The goal of this section is to prove the following result (see Figure
\fig, where the numbers within parentheses refer to the assertions of
the theorem):

\nom\pal
\thm{Theorem \sta} Let $\Gamma $ be a reduced labelled  graph
representing $G$. Suppose $\Out(G)$ does not contain $F_2$. Then:
\roster
\item A slid edge is either a $(2,2)$-segment or a 
$(1,q)$-loop. 
\item  Slid  edges are disjoint. 
\item A pseudo-ascending loop is  a $(p,\pm
p)$-loop or a
$(1,q)$-loop.
\item If $v$ is the basepoint of a $(1,q)$-loop with $|q|\ge2$, then no
other label at
$v$ divides a power of $q$. 
\item
If $v$ is the basepoint of a $(1,q)$-loop with $|q|\ge2$, and  $r,s$ are two
labels at $v$ not carried by the loop, then
$s$ does not divide any $rq^n$. 
\item Let $vw$ be a $(2,2)$-segment. Let $r$ be a label at $v$, and $s$ a
label at $w$ (other than those carried by $vw$). If $r|s$ and $s$ is even, then
$r=s$ and the labels are carried by the same non-oriented edge.  

\endroster
\fthm

We prove Theorem \pal{} in several steps. Our two main tools will be
slide moves and Lemma \tw.

\subhead Slid segments are   $(2,2)$-segments   \endsubhead

 Let $f$ be a slid  
$(p,q)$-segment with $p,q\ge2$.  Consider    the graph of groups $\Theta $
obtained by collapsing $f$. It has a vertex group $H_v$
isomorphic to $\langle a,b\mid a^p=b^q\rangle$. Let
$e$ be an edge of
$\Gamma $ that may slide across $f$, viewed as an edge of $\Theta $.
Its group $G_e$ is generated by a power of $a^p$, so is central in $H_v$. 
By Lemma \tw, the group $\T(\Theta )\inc\Out(G)$ maps onto $J=
Z_{H_v}(G_e)/\langle Z( H_v),Z(G_e)\rangle=\langle a,b\mid
a^p=b^q=1\rangle$. If
$(p,q)\neq(2,2)$, the group $J$ contains
$F_2$, so $\Out(G)$ contains $F_2$. 

\nom\tee
\example {Remark \sta} For future reference, note that $\T(\Theta )$ contains
$\T(\Gamma )$ by Remark \te, and furthermore the image of $\T(\Gamma )$ in
$J$ is finite. To see this, recall that $\T(\Gamma )$ is generated by the groups
$Z_{G_{o(f')}}(G_{f'})$ with $f'\neq f,\bar f$ (Remark \te). All these groups have
trivial image in
$J$, except
$Z_{G_{o(e)}}(G_e)$ whose image is finite.
\endexample

\subhead Slid loops are     
$(1,q)$-loops \endsubhead

We assume that $e$ slides across a   $(p,q)$-loop $f$ with $2\le
p\le|q|$, and we show that $\Out(G)$ contains $F_2$. If
$|q|\ge 3$, we  create a slid $(p,q)$-segment  by performing an expansion
move (replacing the loop by a $(p,q)$-segment and a $(1,1)$-segment), 
and we apply the previous argument (which is valid even if the labelled graph
is not reduced). If
$f$ is a
$(2,\pm2)$-loop, we collapse it and we apply Lemma \tw. We now have
$H_v=\langle a,t\mid ta^2t\mi=a^{\pm2}\rangle$, and $H_e$ is
generated by  a power of $a^2$. The quotient  $Z_{H_v}(G_e)/\langle Z(
H_v),Z(G_e)\rangle$ is isomorphic to $\Z*\Z/2\Z$, so contains $F_2$. 

\subhead  Pseudo-ascending loops are $(p,\pm
p)$-loops or  
$(1,q)$-loops  \endsubhead

This amounts 
to showing that $\Gamma $ cannot contain an $(r, rs)$-loop with $r\ge2$
and $|s|\ge 2$. If it does, write
$ G=G_1*_{\langle a\rangle}G_2$, with $G_1=\langle a,t\mid
ta^rt\mi=a^{rs}\rangle$. By [\Collev], there exist two automorphisms of
$G_1$ fixing $a$ and generating a free subgroup of rank 2 in
$\Out(G_1)$ (in the notation of [\Collev],   set  $\alpha ^r=\varphi
_0=\tau =1$ to see that the subgroup of
$\Out(G_1)$ generated by $\alpha $ and $\gamma _2$ maps
  onto $\Z/r\Z*\Z/s^2\Z$). Extend the automorphisms by the identity
on $G_2$ and check that they generate $F_2\inc\Out(G)$. 

Here is another argument, valid when $r$ or $|s|$ is bigger than 2:
perform an expansion and a slide to obtain a graph with a slid 
$(r,s)$-segment.

\subhead Slid edges are disjoint \endsubhead

We argue by way of contradiction. 
There are several cases to consider (they are pictured from right to left on
Figure 8). First suppose that
$v$ belongs  to a slid
$(2,2)$-segment $f$ and a $(1,q)$-loop $f'$.  If $q$ is even, one may  slide $f'$ 
across $f$ and then collapse $f'$. This creates a $(2,2q)$-loop, a contradiction.
If $q$ is odd, some other edge
  may slide across $f$.  Sliding   $f$ around $f'$
makes $f$ a slid $(2,2q)$-segment,  a contradiction if $|q|>1$. If $q=\pm1$,
collapse both 
$f$ and
$f'$ and apply Lemma
\tw.

Now suppose  $v$ belongs to a  $(1,q)$-loop $f$ and a  $(1,r)$-loop $f'$.
If $|q|\ge2$, sliding $f'$ around $f$   makes it  a $(q, r)$-loop, and then sliding
 $\bar f'$ twice makes it a
$(q,q^2r)$-loop.  
If 
$|q|=|r|=1$, we may write $G$ as an amalgam $G_1*_{\langle
a\rangle}G_2$, with $G_2=\langle a,t,t'\mid
tat\mi=a^{\pm1},t'at'{}\mi=a^{\pm1}\rangle$. It is easy to embed $F_2$
into $\Out(G)$ by using automorphisms of $G_2$ fixing $a$. 

Finally, suppose that $(2,2)$-segments $f$ and $f'$ have a vertex $v$ in
common. We may assume that their other endpoints are distinct, as
otherwise sliding $f$ across $f'$ would create a  
slid  $(2,\pm2)$-loop.   

The fundamental group of the subgraph of groups $f\cup f'$ is 
$J=\langle a,b,c\mid a^2=b^2=c^2\rangle $. Consider the following
automorphisms
$\alpha ,\beta $ of $J$: $\alpha $ fixes $a$ and $b$ and conjugates $c$ by
$ba$, while
$\beta
$ fixes $b,c$ and conjugates $a$ by $bc$. They extend to automorphisms
of $G$ (they are  twists in the graph of groups obtained from $\Gamma
$ by collapsing  
$f'$ and
$f$ respectively). We claim that $\alpha ,\beta $ generate a free
nonabelian subgroup of $\Out(J)$ (hence also of
$\Out(G)$ because $J$ is its own normalizer). 

Indeed, consider $\ov J\simeq
\Z/2\Z*\Z/2\Z*\Z/2\Z$  obtained by adding the relation $a^2=1$. Let
$J^+\inc \ov J$ be the subgroup of index $2$   consisting of elements
of even length. It is free with basis
$ \{\bar a\bar b,\bar b\bar c\}$. With respect to this basis,  $\alpha $
acts on  the abelianization of
$J^+$ as the matrix  $\pmatrix 1 & 2 \\ 0 & 1 
\endpmatrix$,   $\beta $ acts as $\pmatrix 1 & 0 \\ 2 & 1 
\endpmatrix$, and inner automorphisms of $\ov J$ act as $\pm Id$. It
follows that there is no nontrivial relation between $\alpha $ and
$\beta $ in $\Out(J)$.

\subhead Labels near a $(1,q)$-loop  \endsubhead

Let $v$ be the basepoint of a $(1,q)$-loop. We already know that no
other label $r$ at $v$ equals 1 or divides $q$ (it would be carried by a
slid edge).  Suppose that  
$r$ divides some $q^n$. Let $\ell$ be a prime divisor of $r$.  Expand $v$
so as to create a
$(1,\ell)$-segment, and collapse the $(1,\frac q \ell)$-edge (see Figure
\fig).
The new labelled graph is isomorphic to $\Gamma
$, except that $r$ has been divided by $\ell$   and other indices at $v$
have been multiplied by  $\frac q \ell$. Repeat this operation until
$r$ divides $q$, a case already ruled out.

\midinsert
\centerline 
{\includegraphics[scale=.45] 
{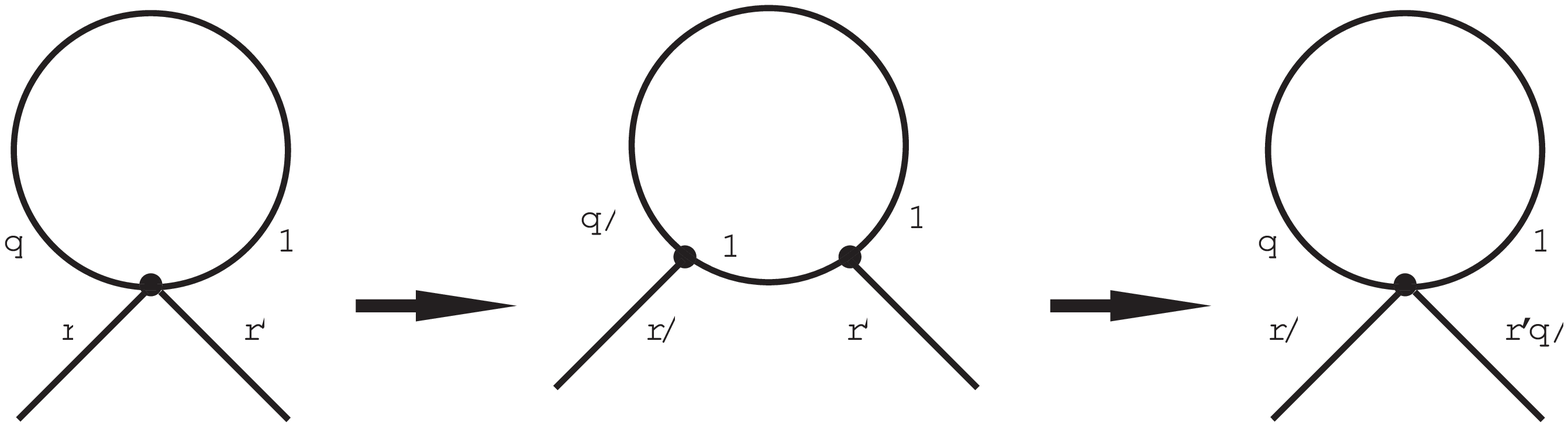}}
\captionwidth{220pt}
\botcaption
 {Figure \the\figno}
{Replacing $r$ by $r/\ell$. }
\endcaption
\endinsert

Now let  $r,s$ be labels at $v$ carried by edges $e,f$, with   $s|rq^n$.    If
$f\neq\ov e$, we can make $e$ a slid edge by sliding it 
$n$ times across the loop.   If $f=\ov e$, we can create an  $(s,
rq^{n+1})$-loop, contradicting (3).

\subhead Labels near a $(2,2)$-segment \endsubhead

Let $vw$ be a $(2,2)$-segment. We already know that all other  labels at
$v$ and
$w$ are bigger than 2 in absolute value. Furthermore, if $r,s$ are labels
at the same vertex ($v$ or $w$), and
$r|s$, then $r=\pm s$ and they are carried by a loop. Assertion (5) of
the theorem follows:  since $s$ is even, one can perform a slide across
$vw$ so that
$s$ becomes a label at $v$. 

\head  \sect { Groups with  $\Out(G)\not\supset F_2$ }\endhead

In this section, we complete the proof of Theorem \six{} by  showing:

\nom\main
\thm{Theorem \sta} If a non-elementary GBS group $G$ is represented by a
reduced labelled graph
$\Gamma $ satisfying the conclusions of Theorem  \pal, then 
$\Out(G)$ is virtually nilpotent of class at most 2.
\fthm

The theorem is true if $G$ is a solvable Baumslag-Solitar group [\Coll], so we
rule out this case.  
As we wish to study the whole automorphism group of $G$, it is important in
this proof to think of
$\Gamma
$   as a marked graph. As usual, we denote by $T$   its Bass-Serre tree.

In general, there are many graphs representing $G$, and the first step  in the
proof of Theorem \main{} will be to show  that  collapsing the slid edges of
$\Gamma $ yields a marked graph of groups $\Theta $ and a (non GBS)
tree $S$ which are   canonical (they do not depend on $\Gamma $). In the
language of [\Forr], we first show that all possible trees $S$ belong to the
same deformation space, and then using [\Ler] that there is only one reduced
tree in that space.  In particular
$\Out(G)=\Out^S(G)$, and we shall conclude the proof of Theorem \main{} by
showing that
$\Out^S(G)$ is virtually nilpotent.

We know  that the slid edges of $\Gamma $ are disjoint, and
are either
$(2,2)$-segments or $(1,q)$-loops. Define $\Theta $ and $S$ by collapsing
them.

Consider edge and vertex groups of $\Theta $.
Edge groups are cyclic. Non-cyclic vertex groups are
Klein bottle groups
$\langle a,b\mid a^2=b^2\rangle$ arising from collapsed $(2,2)$-segments,
and solvable Baumslag-Solitar groups
$BS(1,q)= \langle a,t\mid tat\mi=a^q\rangle$ arising from collapsed
$(1,q)$-loops. Note two special cases:
$\Z^2$ if $q=1$, a Klein bottle group if $q=-1$. 

It is useful to think of $ BS(1,q)$, for $|q|\ge2$, as the subgroup of the affine
group of $\R$
  generated by $a:x\mapsto x+1$ and
$t:x\mapsto qx$. It consists of all maps of the form $x\mapsto q^\alpha
x+\beta $ with $\alpha \in\Z$ and $\beta \in\Z[\frac1{|q|}]$. One deduces,
for instance, that powers
$a^i, a^j$  are conjugate  if and only if
$\frac ij$ is a power of $q$. The element $a^r$ has an $s$-root if and only if $s$
divides some $rq^n$; the root is then unique.

\nom\inva
\thm{Lemma \sta} The set of vertex stabilizers of $S$ does not depend  on the
marked graph  
$\Gamma $. In particular, it is $\Out(G)$-invariant. 
\fthm

\demo{Proof} By Corollary \ellind, it suffices to show that   the vertex
stabilizers of
$S$ are determined by the   elliptic subgroups of $T$.  This will be done
using the tree $T$, but the description of vertex stabilizers will involve
only the algebraic structure of the set of elliptic subgroups. 

Let $v$ be a vertex
of
$T$, and
$\ov v$ its projection in
$\Gamma $. We denote by $\pi :T\to S$ the
collapse map.  We want to understand the stabilizer
$H_v$ of $\pi (v)$ in $S$. It is also the stabilizer of the subtree $\pi \mi(\pi
(v))\inc T$.  

If no  strict ascending loop is attached at $\ov v$, the
stabilizer
$G_v$ of
$v$ is a maximal elliptic subgroup of $T$. Conversely every  maximal
elliptic subgroup arises in this way.   
If a $(1,q)$-loop with $|q|\ge2$ is
attached at $\ov v$, there is no maximal elliptic
subgroup containing $G_v$; if $T'$ is another GBS tree,
then
$G_v$ is not necessarily a vertex stabilizer of $T'$ (see [\For]). 

We first determine  the group
$H_v$   containing a   maximal elliptic subgroup $G_v$.

Consider the normalizer $N(G_v)$. If it is $\Z^2$ or a
Klein bottle group, then $\ov v$ bounds a $(1,\pm1)$-loop and
$H_v=N(G_v)$. Otherwise, $N(G_v)=G_v$. Let $a$ be a generator of
$G_v$.  Then 
$H_v$ is   the centralizer $Z(a^2)$ if $\ov v$ bounds a slid
$(2,2)$-edge,  
$G_v$ if not. 

To decide which (in terms of $G_v$ only), first observe that $\ov
v$ bounds a  
$(2,2)$-edge (slid or not) if and only if $Z(a^2)$ is a Klein bottle group.
Assuming it is, consider (as in [\GLL]) the set of groups of
the form
$Z(a^2)\cap K$, where
$K$ is an elliptic subgroup of $T$ not contained in $Z(a^2)$. It is easy to see
that the edge is slid if and only if some maximal element of
this set (ordered by inclusion) is contained in $\langle a^2\rangle$.

Now suppose that a $(1,q)$-loop with $|q|\ge2$ is
attached at $\ov v$.
Then   $H_v= \langle a,t\mid tat\mi=a^q\rangle$,
where $a$ is a generator of $G_v$.  

Condition (4) of Theorem \pal{} implies that
{\it $\pi (v)$ is the only point of $S$ fixed by $a^\ell$ if $\ell$ divides a power
of
$q$\/}: stabilizers of edges adjacent to $\pi (v)$ in $S$ are conjugate in
$H_v$ to
$\langle a^r\rangle$, where $r$ does not divide any $q^n$; since   $a^\ell$
is not   conjugate (in $H_v$) to a power of such an $a^r$, it cannot fix an edge. 

As in [\For],   say that an
elliptic subgroup $K$ of $T$ is {\it vertical\/} if any elliptic subgroup $K'$
containing $K$ is contained in  a conjugate of
$K$. For the action of $H_v$ on $T$, {\it  the subgroup  $\langle  a^\ell\rangle$
is vertical if and only if 
$\ell$ divides a power of $q$\/}. We show that the   same result holds
for the action of $G$ on $T$. 

Suppose that $\ell$
divides a power of $q$. If   $K'\supset \langle  a^\ell\rangle$ is elliptic (in
$T$, hence also in $S$), then $\pi (v)$ is the only point of $S$ fixed by
$K'$, so $K'$ fixes a point $w\in \pi \mi(\pi (v))\inc T$. The
stabilizer of $w$ is conjugate to $\langle  a \rangle$ in $H_v$, and
$\langle  a \rangle$ is contained in  a conjugate of $\langle  a^\ell \rangle$
because $\ell $ is a power of $q$. Thus $K'$ is contained in a conjugate of
$\langle  a^\ell\rangle$, and $\langle 
a^\ell\rangle$ is vertical as required. Conversely, if
$\langle  a^\ell\rangle$ is vertical, then it   contains   a conjugate  
$g\langle  a \rangle g\mi$ with $g\in G$. Since $\pi (v)$ is the only point of $S$
fixed by
$a $, we have $g\in H_v$ and we deduce that $\ell$ divides a power of $q$.

We now conclude the proof, by characterizing the vertex stabilizer  $H_v$ of
$S$ containing a vertical subgroup
$K\inc G$   which is not maximal elliptic. We know that
$K$ is generated  by $ a^\ell $, where $a$ generates a vertex stabilizer $G_v$,
a
$(1,q)$-loop with $|q|\ge2$ is attached to $\ov v$, and $\ell$ divides a
power of
$q$.   In particular,
$K$ fixes a unique
point $\pi (v)\in  S$. 

The stabilizer $H_v$ of $\pi (v)$ is isomorphic to $BS(1,q)$. The set of
elements of
$H_v$ which are elliptic in $T$ is an abelian subgroup (isomorphic to
$\Z[\frac1{|q|}]$), so
$gKg\mi$ commutes with
$K$ if
$g\in H_v$.
Conversely, if $gKg\mi$ commutes with
$K$, then $g\in H_v$ because $\pi (v)$ is the only fixed point of $K$ in
$S$. We
may now  characterize $H_v$ (independently of
$T$) as the set of
$g\in G$ such that
$gKg\mi$ commutes with
$K$.
\cqfd\enddemo
 
\nom\indep
\thm{Lemma \sta} The $G$-tree $S$ does not depend on $\Gamma $. 
In particular, 
$\Out^S(G)=\Out(G)$.
\fthm

\demo{Proof}  We apply the main result of [\Ler].
Since $S$ is reduced (no inclusion $G_e\hookrightarrow H_v$ is onto), it
suffices to check that the following holds. Let $e$ and
$f$ be oriented edges of
$S$ with the same origin   such that $G_f\inc G_e$; if $e,f$ do not
belong to the same
$G$-orbit,   then  $e,\ov f$ are in the same orbit and $G_e=G_f$. 

Let $v,w$ be the origins of $e$ and $f$ in $T$.   Let $  e_0$, $  f_0$ be the
projections 
  in $\Gamma $, and $r,s$ the corresponding labels. We distinguish several
cases.

First assume $\ov v=\ov w$ in $\Gamma $.  
 If no  collapsing takes place at
$\ov v$, or if $\ov v$ bounds a $(1,\pm 1)$-loop, then $G_f\inc G_e$ implies
that $r$ divides $s$. This is possible only if $e_0$ and $f_0$ are opposite
edges forming a  
$(p,\pm p)$-loop.  

Now suppose that $\ov v$ bounds a $(1,q)$-loop with $|q|\ge2$. Write
the corresponding vertex stabilizer $H_v$ of $S $ as 
$ \langle a,t\mid tat\mi=a^q\rangle$. Then $\langle  a^s\rangle$ is conjugate
in $H_v$ to a subgroup of $\langle a^r
\rangle$, so there exists $n$ such that $\frac s{nr}$ is a power of $q$.
This contradicts Assertion (5) of Theorem \pal, so this case cannot
occur.

Finally, suppose that  $\ov v\ov w$ is a slid $(2,2)$-edge.
The stabilizer of $\pi (v)$ in $S$  is then $ H_v=\langle
a,b
\mid a^2=b^2\rangle$. The   subgroups $G_e$ and $G_f$ of $H_v$ are generated
by conjugates of  powers of $a$ or $b$.  Distinct
powers of
$a$ (resp\. $b$) are not conjugate in $H_v$, while $a^i$  is conjugate to  $b^j$
only when $i$ and $j$ are equal and even.  In particular, $r$ 
 divides $s$. If $e_0$ and $f_0$ have the same
origin ($\ov v$ or $\ov w$),   we conclude as in the first case. If not, then $s$ 
must be even and we use Assertion (6) of Theorem \pal. 
\cqfd\enddemo

We may now study $\Out(G)=\Out^S(G)$ using the results of [\LGD]
recalled in Section 3. In particular, $\rho $ has a restriction 
$\rho _1:\Out_1^S (G)\to\prod_{u\in W}\Out(H_u)$ with $\Out_1^S (G)$ of
finite index in $\Out(G)$ and $\ker\rho _1= \T(S)$ (we denote by $W$ the
vertex set of
$\Theta $, and by $H_u$ the vertex group of $u\in W$).   

 We first show that $\T(S)$ is virtually abelian. It is generated by centralizers
of edge groups in vertex groups $H_u$. If $H_u$ is $\Z$ or $\Z^2$, the
centralizer is of course $H_u$. If
$H_u=BS(1,q)$ with
$|q|\ge2$, the centralizer is an infinitely generated abelian
group isomorphic to $\Z[\frac1{|q|}]$. If
$H_u$ is a Klein bottle group, the centralizer is $\Z^2$  if $u$ comes
from   a $(1,-1)$-loop, $\Z$ if $u$  comes from a slid $(2,2)$-segment   and
the edge group is not central, the whole of $H_u$ if the edge group is
central. Since a Klein bottle group is virtually abelian, so is $\T(S)$. 

\nom\fingen
\example{Remark \sta} Relations in the presentation of $\T(S)$ come from 
centers of edge and vertex groups of
$\Theta
$. Since these centers are $\{1\}$, $\Z$, or $\Z^2$, the group
$\T(S)$ is finitely generated if and only if $\Gamma $ contains no strict
ascending loop.
\endexample

Now fix a vertex
$u$ of $\Theta $, and define    $P_u\inc \Out(H_u)$  by projecting  the
image of
$\rho _1$. 
If
$H_u$ is
$\Z$ or a Klein bottle group, $P_u$ is finite because $\Out(H_u)$ is
finite. We claim that $P_u$ is finite also when $H_u$ is $BS(1,q)$ with
$|q|\ge2$.

Write $
H_u=\langle a,t\mid tat\mi=a^q\rangle$. The vertex $u$ is obtained by
collapsing a $(1,q)$-loop $f_u$ of $\Gamma $. Denote its basepoint by
$v$.
Since $G$ is assumed not to be solvable, we may consider an edge group
$\langle a^r\rangle$, where $r$ is a label near $v$ not carried by $f_u$. Its
image by an automorphism
$\alpha 
\in P_u$ is also an edge group, so $\alpha (\langle a^r\rangle)$  is conjugate to
$\langle a^s\rangle$ for some label $s$ (possibly equal to $r$). But $a^r$ has
an $s$-th root only if $s$ divides some $rq^n$, so  
$r=s$ by Assertion (5) of Theorem \pal. By uniqueness of roots,  $\alpha $
maps
$a$ to a conjugate of
$a^{\pm1}$. Only finitely many outer automorphisms of $H_u$ have this
property [\Coll], so $P_u$ is indeed finite. 

The group $P_u$ is infinite only when $u$ comes from collapsing a $(1,1)$-loop
${f_u}$. In this case, $H_u=\langle a,t\mid tat\mi=a \rangle$. As above, $a$
must be mapped to a conjugate of $a^{\pm1}$, so $P_u$ contains with index at
most 2 the group generated by the automorphism $D_{f_u}$  fixing $a$ and
mapping $t$ to $at$. We view $D_{f_u}$ as an    automorphism of $G$ (extend
it by the identity). It is a twist relative to $\Gamma $, but not to $\Theta $
(Remark \te{} does not apply here, as $f_u$ is not a segment; in general, none
of the groups
$\T(S )$, $\T(T )$ contains the other).

This analysis shows that {\it
$\T(S)$ and the automorphisms
$D_{f_u}$ associated to  
$(1,1)$-loops of $\Gamma $  generate a   finite index subgroup   of
$\Out(G)$\/}. 
We replace $\T(S)$ by an abelian subgroup $\T_0(S)$ of finite index,
and
we complete the proof of Theorem \main{} by  showing that the subgroup
generated by
$\T_0(S)$ and the automorphisms $D_{f_u}$  is virtually nilpotent of
class
$\le2$: every commutator is central. 

Non-commutativity only comes from the fact that $D_{f_u}$ may fail to
commute with $D(z)$, when
$z\in H_u$ and
$D(z)$ is a twist of $\Theta $ around an edge
$e$ with origin $u$. Write $H_u=\langle a,t\mid tat\mi=a \rangle$. 
 The group
$G_e$ is generated by a power of
$a$.

Recall that $D_{f_u}$ fixes $a$ and maps $t$ to $at$. In particular, $D(z)$
commutes with $D_{f_u}$ if $z$ is a power of $a$ (both automorphisms belong
to
$\T(T)$). The interesting case is when $z=t$ (geometrically, $u$ carries a
2-torus $T^2$, $e$ carries an annulus attached to a meridian of $T^2$,
$D_{f_u}$ is a Dehn twist in $T^2$ around a meridian, and
$D(t)$ drags the annulus around    $T^2$ along the longitude). But
conjugating $D(t)$ by
$D_{f_u}$ gives $D(ta)$, so the commutator of $D(t)$ and $D_{f_u}$ is $D(a)$, a
central element. This easily implies that every commutator is
central, completing the proof of Theorem
\main.

\head \sect Further results\endhead

\subhead  Nilpotent vs abelian\endsubhead

\thm{Corollary \sta} If $G$ is represented by a reduced labelled graph with no
$(1,1)$-loop, then $\Out(G)$ contains $F_2$ or is virtually abelian.
\fthm

This follows immediately from the proof of Theorem \main. More generally,
if $\Out(G)$ does not contain $F_2$, it is  
virtually abelian if and only if every commutator
$D(a)$ as in the last paragraph of the proof has finite order. This happens in
particular  when the basepoint of every 
$(1,1)$-loop   has valence 3.
See Remark \foo{} for a more complete discussion. 

\subhead  Finite generation\endsubhead

Here is a general fact:

\nom\finig
\thm{Proposition \sta} Let $\Gamma $ be a labelled graph representing
a GBS group $G$. If $\Gamma $ contains  a strict ascending loop, but $G$
is not a solvable Baumslag-Solitar group, then $\Out(G)$ has an
infinitely generated abelian subgroup. 
\fthm

\demo{Proof}   Collapse the loop and apply Lemma
\tw{} to an edge $e$ with   origin at the collapsed vertex $v$ (there is such an
edge because $G$ is not solvable). We have
$H_v=BS(1,q)$ with $|q|\ge2$,   and $Z_{H_v}(G_e)$ is
infinitely generated abelian (it is isomorphic to $ \Z[\frac1{|q|}]$). The
center of $H_v$ is trivial, and the center of $G_e$ is cyclic. The
subgroup of $\T(\Theta )$ generated by $Z_{H_v}(G_e)$ is   isomorphic
to $ \Z[\frac1{|q|}]$, or is an infinite abelian torsion group. 
\cqfd\enddemo

From the proof of Theorem \main{} we get:

\thm{Corollary \sta} Let $\Gamma $ be a reduced labelled graph
representing a non-solvable GBS group $G$ with   $\Out(G)$  
virtually nilpotent. The group $\Out(G)$ is finitely
generated if and only if $\Gamma $ contains no strict ascending loop. 
\fthm

\demo{Proof} We have seen that $\Out(G)$ is generated by the union of
$\T(S)$ and a finite set, so by virtual nilpotence $\Out(G)$ is finitely
generated if and only if
$\T(S)$ is finitely generated. The corollary now follows from Remark \fingen.
\cqfd
\enddemo

\thm{Corollary \sta} If no label of $\Gamma $ equals 1, then $\Out(G)$ contains
$F_2$ or is finitely generated and virtually abelian. \cqfd
\fthm

In the virtually abelian case, the  torsion-free rank   may be
computed from $k$ and the labels near the $(2,2)$-edges.

\subhead  Algebraic rigidity\endsubhead

\thm{Theorem \sta}  If the GBS group $G$ is not a solvable
Baumslag-Solitar group, the following are equivalent:
\roster
\item  $G$ is algebraically rigid (there is only one reduced GBS  tree).
\item The deformation space $P\D$ is a finite complex.  
\item $\Out(G)$ is virtually $\Z^k$ (with $k$ defined in Proposition \kk). 
\item Let $\Gamma $ be any reduced labelled graph representing   $G$. If $e,f$
are distinct oriented edges of $\Gamma $ with the same origin
$v$, and the label of $f$ divides that of $e$, then either $e=\ov f$ is a
$(p,\pm p)$-loop with $p\ge2$, or $v$ has valence 3 and bounds a
$(1,\pm1)$-loop. 
\endroster 
\fthm

\example{Remarks} $\bullet$ If $G$ is unimodular, $(3)\Leftrightarrow(4)$
follows from  Theorem \esout{}
and
  [\Pett, Corollary 5.14].

$\bullet$  Suppose $|n|\ge2$. Then $BS(1,n)$ is algebraically rigid if and only if
$|n|$ is prime [\Ler], while
$\Out(BS(1,n))$ is virtually $\Z^k$ (i.e\. finite) if and only if $|n|$ is a prime
power [\Coll].
\endexample

\demo{Proof}  The equivalence $(1)\Leftrightarrow(4)$ is in [\Ler], and
$(1)\Rightarrow(2)\Rightarrow (3)$ follows from Sections 3 and 5. We prove 
$(3)\Rightarrow (4)$.

Suppose
that $\Out(G)$ is virtually $\Z^k$  (equivalently, $\T(\Gamma )$ has finite
index in $\Out(G)$). Let 
$e,f$ be adjacent edges with
$\lambda _f|\lambda _e$. By Theorem \pal{} and Proposition \finig, the edge
$f$ must be a slid $(2,2)$-segment or     a $(p,\pm p)$-loop.

 It cannot be a
segment because of Remark \tee: after collapsing $f$, the group $\T(\Theta
)\inc\Out(G)$ would map onto the infinite dihedral group $J=\langle a,b\mid
a^2=b^2=1\rangle$ with the image of
$\T(\Gamma )$ finite, a contradiction. 
To prove (4), there remains to show that the basepoint   of any $(1,\pm
1)$-loop  has valence 3.

Let $f$ be a $(1,\varepsilon )$-loop, with $\varepsilon =\pm1$, let $v$ be its
basepoint, let 
$e_1,\dots, e_n$ be the oriented edges with origin $v$  (other than $f,\bar f$).
We must show $n=1$. 

First consider the subgroup
$\T_0$ of
$\T=\T(\Gamma )$ generated by the  twists $D_{e_i}$ and the twists around
edges with origin other than $v$. The group $\T$ is generated by $\T_0$ and
the twists
$D_f$, $D_{\bar f}$. The only relations involving $D_f$, $D_{\bar f}$ are
$D_f+\varepsilon  D_{\bar f}=0$ (edge relation) and $D_f+D_{\bar
f}+\sum_iD_{e_i}=0$ (vertex relation). It follows that $\T_0$ has index at
most 2 in $\T$ if $\varepsilon =-1$, that $\T$ is the direct sum of $\T_0$ and
the infinite cyclic group generated by $D_f$ if $\varepsilon =1$. 

Now let  $\Theta $ be the graph of groups obtained by collapsing $f$, and  
consider  $\T'=\T(\Theta )$.   We will see that it  is generated by
$\T_0$ together with extra twists $D'_i$ around the edges
$e_i$ (note that edge group centralizers are bigger in $\Theta $ than in
$\Gamma $). To describe $D'_i$ precisely, we distinguish two cases.

If $\varepsilon =-1$, the vertex group $H_v$ of $v$ in $\Theta $ is a Klein
bottle group $\langle a,t\mid tat\mi=a\mi\rangle$. Its center is $\langle
t^2\rangle$. The edge groups
$G_{e_i}$ are generated by  powers of $a$. Their centralizer in $H_v$ is the
free abelian group generated by $a$ and $t^2$. In this case, $D'_i$ is the twist
by
$t^2$ around $e_i$. The only relation involving $D'_i$ is the vertex relation
$\sum_iD'_i=0$, so $\T'$ is the direct sum of $\T_0$ and $\Z^{n-1}$. As $\T_0$
has index at most 2 in $\T$, we must have $n=1$ since  $\T$ has
finite index in
$\Out(G)$.

If $\varepsilon =1$, we still have $\T'=\T_0\oplus\Z^{n-1}$ (the vertex group
$H_v$ is $\Z^2=\langle a, t\rangle$, and
$D'_i$ is the twist by $t$).  There is a natural homomorphism from
$\langle\T,\T'\rangle$  to
$\Out(H_v)$, given by the action on
$H_v$.
The kernel contains $\T'$, but its intersection with $\T$ is $\T_0$ (as $D_f$
acts on $H_v$ by fixing $a$ and mapping $t $ to $at$). Since
$\T$ has finite index in $\Out(G)$, hence in $\langle\T,\T'\rangle$,  we deduce
that
$\T_0=\T\cap \T'$ has finite index in
$\T'$, so $n=1$. 
\cqfd\enddemo

Combining with Proposition \kk, we obtain:

\nom\outfini
\thm{Corollary \sta} Let $\Gamma $ be a reduced labelled graph
representing a group  $G$. The group $\Out(G)$ is finite
if and only if one of the following holds:  
\roster
\item
$\Gamma
$ is a   tree with no divisibility relation;
\item
$\Gamma $ is a   graph with first Betti number 1, there is  no divisibility
relation,  and $G$ has nontrivial modulus;
\item
$\Gamma
$ is obtained from a tree with no divisibility relation by attaching one $(k,-
k)$-loop. If
$k\ge2$, no other index at the attaching point is a multiple or a divisor of
$k$. If
$k=1$, the loop is attached at a terminal vertex.
\cqfd
\endroster
\fthm

\subhead On the isomorphism problem
\endsubhead

Given a labelled graph $\Gamma $, it is easy to decide algorithmically
whether the associated GBS group $G$
 is elementary, solvable, unimodular. By Theorems \pal{} and \main, we may
decide whether $\Out(G)$ contains $F_2$ or is virtually nilpotent.

The isomorphism problem for GBS groups is the problem of deciding whether
two (reduced) labelled graphs represent isomorphic groups or not. It is
solved for rigid groups (obviously),   for groups with no nontrivial
integral modulus [\For], and for $2$-generated groups [\Lrang].

\thm{Theorem \sta} 
The isomorphism problem is solvable for GBS groups such
that $\Out(G)$ does not contain a non-abelian free group. 
\fthm

\demo{Proof} Let $\Gamma $ be a reduced labelled graph representing
$G$. We assume that $\Out(G)$ does not contain $F_2$, so $\Gamma $ satisfies
all six conditions of Theorem \pal.

We describe three ways of producing new labelled graphs representing  $G$
(besides admissible sign changes). The first one changes the graph, the other
two only change the labels.

\roster
\item 
  Sliding an edge across a $(2,2)$-segment: it  changes the attaching
point of an edge (carrying an even label).
\item If $v$ is the basepoint of a $(1,q)$-loop, one may multiply or divide
by $q$ some other  label near $v$, by sliding the corresponding edge around
the loop. 

\item If $v$ is the basepoint of a $(1,q)$-loop, one may multiply or divide
all other labels near $v$ by any number $p$ dividing $q$, by performing an
expansion at $v$ followed by   a collapse (see Figure
\fig; this is called an induction move in
[\Ler]). 
\endroster

\midinsert
\centerline 
{\includegraphics[scale=.45]
{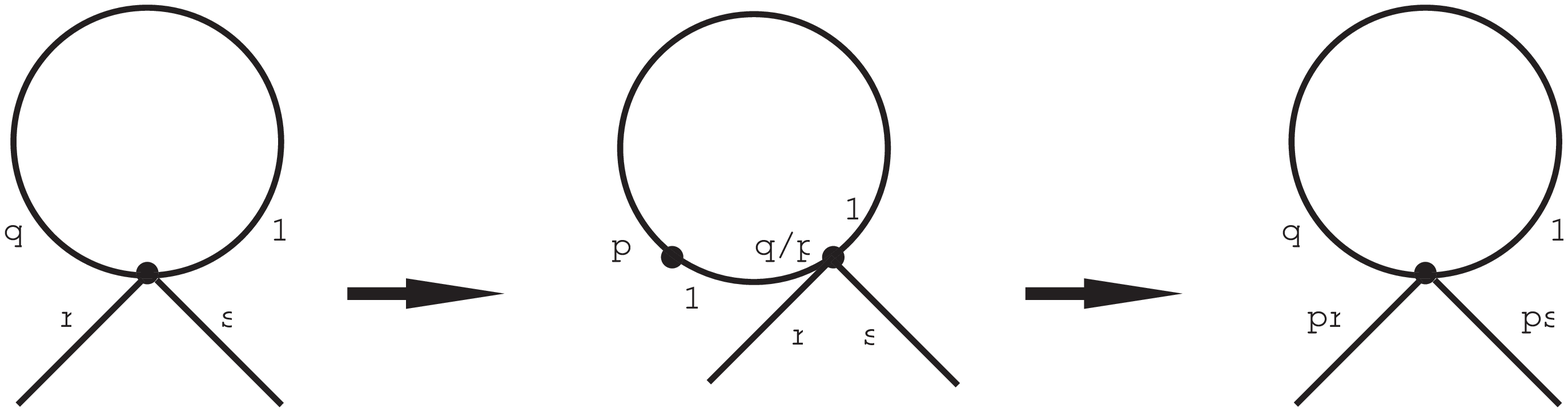}}
\captionwidth{220pt}
\botcaption
 {Figure \the\figno}
{Induction move. }
\endcaption
\endinsert

Consider the set $\G$ consisting of all labelled graphs which may be obtained
from
$\Gamma $ by combining these moves. They are reduced by condition (4) of
Theorem \pal, and it is easy to decide whether a given labelled graph
$\Gamma '$ belongs to $\G$. We now
complete the proof by showing that $\G$ contains all reduced
graphs   representing $G$. 

As above, consider the graph of groups $\Theta $  and the
Bass-Serre tree $S $ obtained by collapsing the slid edges of $\Gamma $.
We have seen   that $S$ does not depend on the graph
$\Gamma
$ used to construct it (Lemma \indep).  It thus suffices to show that the
various ways of blowing up
$S$ into a GBS-tree differ by the moves mentioned above.

First consider a vertex group $H_v$ obtained by collapsing a $(2,2)$-edge. It
is a Klein bottle group $ \langle a,b\mid a^2=b^2\rangle$. The generator
of any adjacent edge group is conjugate to some $a^i$ or $b^j$. Distinct
powers of
$a$ (resp\. $b$) are not conjugate in $H_v$, while $a^i$  is conjugate to  $b^j$
only when $i$ and $j$ are equal and even. It follows that all ways of blowing up
$v$ differ by slides across the $(2,2)$-edge. 

Now consider a vertex group $H_v= \langle a,t\mid tat\mi=a^q\rangle$
obtained by collapsing a $(1,q)$-loop of $\Gamma $. Let $T'$ be another
GBS-tree, associated to a labelled graph $\Gamma '$. As in the proof of
Lemma \inva, we cannot say that $a$ generates a vertex group of $T'$, only
that some $a^i$, with $i$ dividing a power of $q$, does. 

Suppose for a moment
$i=1$. The generator
of any   edge group adjacent to $v$   is conjugate to a
power of
$a$. As
$a^m$, $a^n$ are conjugate in $H_v$ only if $\frac mn$ is a power of $q$, the
labelled graphs $\Gamma ,\Gamma '$ differ by moves of type (2) near the
$(1,q)$-loop. 

If $i\neq1$, it is a product of   divisors of $q$, and $\Gamma
,\Gamma '$ differ by moves of type (2) and (3).
\cqfd\enddemo

\example{Remark} The same technique may be used  when
$G$ is represented by a graph $\Gamma $ satisfying the six conditions of
Theorem \pal, but with arbitrary slid segments allowed in condition (1) (not
only
$(2,2)$-segments). Condition (6) must then be rephrased as follows: Let $vw$
be a $(p,q)$-edge. Let $r$ be a label at $v$, and $s$ a label at $w$. If $q|s$ and
$qr|sp$, then $qr=sp$ and the labels are
carried by the same non-oriented edge.  
\endexample

\bigskip
\Refs 
\widestnumber\no{99} 
\refno=0

\bref \by
 H. Bass, R.  Kulkarni\paper
Uniform tree lattices \jour
J. Am. Math. Soc. \vol3\pages843--902 \yr1990
\endref

\bref\by G. Baumslag, D. Solitar\paper Some two-generator, one-relator
non-Hopfian groups\jour Bull. Amer. Math. Soc. \vol68\yr1962\pages
199--201\endref 

\bref\by M. Clay \paper Contractibility of deformation spaces 
\jour Alg. Geom. Topol. (to appear) \endref

\bref\by M. Clay \paper A generalization of Culler's theorem \jour
arXiv:math.GR/0502248\endref

\bref  \by D.J. Collins \paper The automorphism towers   of some
one-relator  groups\jour Arch. Math.\yr1983\vol40\pages385--400
\endref

\bref  \by D.J. Collins, F. Levin\paper Automorphisms and hopficity of
certain Baumslag-Solitar groups\jour Arch.
Math.\yr1983\vol40\pages385--400
\endref

\bref  \by M. Culler, K. Vogtmann \paper Moduli of graphs and
automorphisms of  free groups\jour Invent. Math. \vol 84
\yr1986\pages91--119 \endref

\bref\by A. Fel'shtyn  \paper  
The Reidemeister number of any automorphism of a Gromov 
hyperbolic group is infinite 
\jour J. Math. Sci.\vol119\pages117-123\yr2004
\endref

\bref\by A. Fel'shtyn, D.L. Gon\c calves \paper Twisted conjugacy
classes of automorphisms of Baumslag-Solitar groups
\jour arXiv:math.GR/0405590\endref

\bref   \by M. Forester\paper Deformation and rigidity of
simplicial group actions on trees\jour Geom.
\& Topol.\vol 6\yr2002\pages 219--267\endref

\bref   \by M. Forester\paper On uniqueness of JSJ
decompositions of finitely generated groups\jour Comm. Math.
Helv. \vol78\yr2003\pages 740--751 \endref

\bref  \by M. Forester\paper Splittings of generalized
Baumslag-Solitar  groups \jour arXiv:math.GR/0502060
\endref

\bref  \by N.D. Gilbert, J. Howie, V. Metaftsis, E. Raptis\paper Tree
actions of automorphism groups\jour J. Group
Theory\yr2000\vol3\pages213--223
\endref

\bref   \by V. Guirardel\paper A very short proof of Forester's
rigidity result\jour Geom. \& Topol.\vol 7\yr2003\pages
321--328\endref

\bref \by V. Guirardel, G. Levitt \paper The outer space   of a free
product \jour arXiv:math.GR/0501288\endref

\bref \by V. Guirardel, G. Levitt \jour in preparation\endref

\bref\by P.H. Kropholler\paper A note on centrality in $3$-manifold
groups\jour Math. Proc. Camb. Phil. Soc.\yr1990\vol107\pages 261--266
\endref

\bref\by P.H. Kropholler\paper Baumslag-Solitar groups and some other
groups of cohomological dimension two\jour Comm. Math. Helv.
 \yr1990\vol65\pages 547--558
\endref

\bref  \by 
    S. Krsti\'c, K. Vogtmann\paper Equivariant outer space and
automorphisms of free-by-finite groups\jour Comm. Math. Helv.
\vol 68 \yr1993\pages 216--262 
\endref

\bref\by G. Levitt \paper Automorphisms of   hyperbolic
groups and graphs of groups \jour Geom.
Dedic. \vol 114\yr2005
\pages 49--70 
\endref 

\bref\by G. Levitt \paper  Characterizing  rigid simplicial actions on
trees
\inbook Geometric methods in group
theory,   Contemp. Math.  \vol372\yr2005\pages27--33
\endref

\bref\by G. Levitt\jour in preparation
\endref

\bref\by G. Levitt, M. Lustig\paper Most automorphisms of a
hyperbolic group have 
very simple dynamics \jour Ann. Sc. ENS\vol33\yr2000
\pages507--517\endref

\bref   \by D. McCullough, A. Miller \paper Symmetric automorphisms
of free products \jour Mem. Amer. Math. Soc.
\vol122\yr1996\endref

\bref
\by D.I. Moldavanskij\paper
Isomorphism of the Baumslag-Solitar groups
\jour Ukr. Mat. Zh. \vol 43\pages1684-1686 \yr1991 
\transl\nofrills \jour Ukr. Math. J. \vol 43\pages 1569-1571 \yr1991 
\endref

\bref \by L. Mosher, M.  Sageev, K. Whyte\paper Quasi-actions on trees,
I. Bounded valence\jour Ann. Math\yr2003\vol158\pages115--164
\endref

\bref  \by M.R. Pettet\paper Virtually free groups with finitely many
outer automorphisms\jour Trans. AMS
 \vol349
\yr1997\pages4565--4587
\endref

\bref  \by M.R. Pettet\paper The automorphism group of a graph
product of groups\jour Comm. Alg.
 \vol27(10)
\yr1999\pages4691--4708
\endref

\bref  \by  E. Raptis, D. Varsos\paper On the automorphism group of the
fundamental group of a graph of polycyclic groups\jour Algebra
Colloquium\yr1997\vol4\pages241--248
\endref

\bref \by A. Rhemtulla, D.  Rolfsen \paper
Local indicability in ordered groups: braids and elementary amenable groups
\jour Proc. Am. Math. Soc. \vol130\pages 2569-2577 \yr2002
\endref

\bref\by K. Whyte\paper The large scale geometry of the higher
Baumslag-Solitar groups\jour Geom. Funct. Anal. \yr
2001\vol11\pages 1327--1343\endref

\endRefs

\address   LMNO, umr cnrs 6139, BP 5186, Universit\'e de Caen,
14032 Caen Cedex, France.\endaddress\email 
levitt\@math.unicaen.fr
\endemail

\enddocument